\documentclass[10pt]{amsart}
\usepackage[]{amsmath}
\usepackage[]{amssymb}
\usepackage{amsfonts}
\usepackage{amsrefs}

\newcommand{\abs}[1]{{\left|#1\right|}}
\newcommand{\set}[1]{{\left\{#1\right\}}}
\newcommand{\st}{\colon}

\newcommand{\RP}{\mathop\mathsf{RP}}
\newcommand{\PR}{\mathop\mathsf{PR}}
\newcommand{\RPPR}{\mathop\mathsf{RPPR}}
\newcommand{\ANY}{\mathop\mathsf{ANY}}
\newcommand{\ORD}{\mathop\mathsf{ORD}}
\newcommand{\ord}{{\mathop{\mathrm{ord}}\nolimits}}
\newcommand{\Li}{{\mathop{\mathrm{Li}}\nolimits}}

\newcommand{\Nfp}{F}
\newcommand{\Mfp}{G}
\newcommand{\Ntc}{T}
\newcommand{\Nha}{C}
\newcommand{\sigzero}{d}
\newcommand{\sigone}{\sigma}

\newtheorem{theorem}{Theorem}[section]
\newtheorem{lemma}[theorem]{Lemma}
\newtheorem{conjecture}[theorem]{Conjecture}
\newtheorem{proposition}[theorem]{Proposition}
\newtheorem{heuristic}[theorem]{Heuristic}
\newtheorem{corollary}[theorem]{Corollary}

{
\theoremstyle{remark}
\newtheorem{remark}[theorem]{Remark}
}

\begin{document}

\title[Small Cycles of the Discrete Logarithm]{Some Heuristics and
Results for Small Cycles of the Discrete Logarithm}

\author{Joshua Holden}

\thanks{The first author would like to thank the Rose-Hulman Institute of
Technology for the special stipend which supported this project 
during the summer of 2002.}

\address{Department of Mathematics,
Rose-Hulman Institute of Technology,
Terre Haute, IN, 47803-3999, USA}

\email{holden@rose-hulman.edu}

\author{Pieter Moree} 

\thanks{The research of the second author was carried out whilst he
was visiting assistant professor at the University of Amsterdam and
supported by Prof.~E.M. Opdam's Pioneer grant of the Netherlands
Organization for Scientific Research (NWO)\@.}

\address{Max-Planck-Institut f\"ur Mathematik, Vivatsgasse 7, D-53111
Bonn, Germany} 

\email{moree@mpim-bonn.mpg.de} 

\date{\today}

\begin{abstract}
Brizolis asked the question: does every prime $p$ have a pair $(g,h)$
such that $h$ is a fixed point for the discrete logarithm with base
$g$?  The first author previously extended this question to ask about
not only fixed points but also two-cycles, and gave heuristics
(building on work of Zhang, Cobeli, Zaharescu, Campbell, and
Pomerance) for estimating the number of such pairs given certain
conditions on $g$ and $h$.  In this paper we extend these heuristics
and prove results for some of them, building again on the
aforementioned work.  We also make some new conjectures and prove some
average versions of the results.
\end{abstract}

\maketitle

\section{Introduction and Statement of the Basic Equations}

Paragraph~F9  of~\cite{UPINT} includes the following problem,
attributed to Brizolis: given a prime $p>3$, is there always a
pair $(g,h)$ such that $g$ is a primitive root of $p$, $1 \leq h
\leq p-1$, and
\begin{equation} \label{fp}
g^{h} \equiv h  \mod{p} \enspace ?
\end{equation}
In other words, is there always a primitive root $g$ such that the
discrete logarithm $\log_{g}$ has a fixed point?  As we shall see,
Zhang~(\cite{Zhang}) not only answered the question for sufficiently
large $p$, but also estimated the number $N(p)$ of pairs $(g,h)$ which
satisfy the equation, have $g$ a primitive root, and also have $h$ a
primitive root which thus must be relatively prime to $p-1$.  This
result seems to have been discovered and proved by Zhang
in~\cite{Zhang} and later, independently, by Cobeli and Zaharescu
in~\cite{CZ}.  Campbell and Pomerance (\cite{CampbellThesis}) made the
value of ``sufficiently large'' small enough that they were able to use
a direct search to affirmatively answer Brizolis' original question.
As in~\cite{Holden02}, we will also consider a number of variations
involving side conditions on $g$ and $h$.

In~\cite{Holden02}, the first author also investigated the two-cycles of
$\log_{g}$, that is the pairs $(g,h)$ such that there is some $a$
between $1$ and $p-1$ such that
\begin{equation}\label{tc}
g^{h} \equiv a \mod{p} \quad \text{and} \quad g^{a} \equiv h \mod{p} .
\end{equation}
As we observed, attacking~\eqref{tc} directly requires the
simultaneous solution of two modular equations, presenting both
computational and theoretical difficulties.  Whenever possible,
therefore, we instead work with the modular equation
\begin{equation} \label{ha}
    h^{h} \equiv a^{a} \mod{p} .
\end{equation}
Given $g$, $h$, and $a$ as in~\eqref{tc}, then \eqref{ha} is clearly
satisfied and the common value is $g^{ah}$ modulo $p$.  Conditions on
$g$ and $h$ in~\eqref{tc} can (sometimes) be translated into
conditions on $h$ and $a$ in~\eqref{ha}.  On the other hand, given a
pair $(h,a)$ which satisfies~\eqref{ha}, we can attempt to solve for
$g$ such that $(g,h)$ satisfies~\eqref{tc} and translate conditions on
$(h,a)$ into conditions on $(g,h)$.  Again, we will investigate using
various side conditions.

Using the same notation as in~\cite{Holden02}, we will refer to an
integer which is a primitive root modulo $p$ as $\PR$ and an integer
which is relatively prime to $p-1$ as $\RP$.  An integer which is both
will be referred to as $\RPPR$ and one which has no restrictions will
be referred to as $\ANY$.  In some instances,  $\bullet$ will be used 
to stand for any one of these four conditions.

All integers will be taken to be between $1$ and $p-1$, inclusive,
unless stated otherwise.  If $N(p)$ is, as above, the number of
solutions to~\eqref{fp} such that $g$ is a primitive root and $h$ is a
primitive root which is relatively prime to $p-1$, then we will say
$N(p)=\Nfp_{g \PR, h \RPPR}(p),$ ($\Nfp$ for ``fixed points'') and
similarly for other equations and conditions.  Likewise the number of
solutions to~\eqref{tc} will be denoted by $\Ntc$ (for ``two-cycles'')
and the number of solutions to~\eqref{ha} will be denoted by $\Nha$
(for ``collisions'').  If $\ord_{p}(g)=\ord_{p}(h)$, we say that $g
\ORD h$.

The first part of this paper focuses on solutions to~\eqref{fp}, with
Section~\ref{indepsec} covering the basic heuristics used and the 
lemmas which can be proven about them.  Section~\ref{fpconjsec} 
presents the conjectures about solutions to~\eqref{fp} which follow 
from the heuristics, and Section~\ref{fpthmsec} proves some new 
theorems which give support to the conjectures.

The middle of the paper deals with solutions to~\eqref{tc} 
and~\eqref{ha}.  Section~\ref{equivsec} examines the relationship 
between solutions of the two equations, while Section~\ref{tcconjsec} 
presents the heuristics used to estimate the number of solutions to 
these two equations and the conjectures that follow from these 
heuristics.

The later sections of the paper deal with average versions of the conjectures 
and results presented in previous sections.  Section~\ref{aveconjsec} 
sets out the lemmas we need and gives average versions of the 
conjectures.  Section~\ref{averessec} gives average versions of the 
results we have proved, where possible, and makes conjectures on the 
others.  Section~\ref{conclsec} discusses further work to be done 
along the lines of this paper.

\section{The ``Independence'' of Order and GCD}
\label{indepsec}

The fundamental observation at the heart of the estimation of
$\Nfp_{g \PR, h \RPPR}(p)$ is that if $h$ is a primitive root
modulo $p$ which is also relatively prime to $p-1$, then there is a
unique primitive root $g$ satisfying~\eqref{fp}, namely $g =
h^{\overline{h}}$ reduced modulo $p$, where $\overline{h}$ denotes the
inverse of $h$ modulo $p-1$ throughout this paper.  Thus to estimate
$N(p)$, we only need to count the number of such $h$; $g$ no longer
has to be considered.  We observe that there are $\phi(p-1)$
possibilities for $h$ which are relatively prime to $p-1$, and we
would expect each of them to be a primitive root with probability
$\phi(p-1)/(p-1)$.  This heuristic uses the assumption that the
condition of being a primitive root is in some sense ``independent''
of the condition of being relatively prime.

\begin{heuristic} \label{rpprheur}
    The condition of $x \RP$ is independent of the condition that $x
    \PR$, in the sense that for all $p$,
    \begin{multline*}
    \frac{\#\set{x \in \set{1, \ldots, p-1} \st x \RPPR}}{p-1} \\
    \approx
    \frac{\#\set{x \in \set{1, \ldots, p-1} \st x \RP}}{p-1} \cdot
    \frac{\#\set{x \in \set{1, \ldots, p-1} \st x \PR}}{p-1}.
    \end{multline*}
\end{heuristic}

That this is essentially
the case was proved in~\cite{Zhang} and in~\cite{CZ}.
We start with the key lemmas of~\cite{CZ}.  Fix a prime $p$.
Let
\[
\mathcal{P}=\mathcal{P}(a,r,N) = \{ a, a+r, \ldots, a + (N-1)r\}
\]
be an arithmetic progression, where $a$, $r$, and $N$ are positive
integers such that $\mathcal{P} \subseteq \{1, \ldots, p\}$.  Let
\[
\mathcal{P}^{\PR}=\{ x \in \mathcal{P} \st x \PR \}
\]
(this is called $\mathcal{P}'$ in~\cite{CZ}),
\[
\mathcal{P}^{\RP}=\{ x \in \mathcal{P} \st x \RP \},
\]
and
\[
\mathcal{P}^{\RPPR}=\{ x \in \mathcal{P} \st x \RPPR \}.
\]
  Finally, for any set of integers $\mathcal{S}$, let
\[ \textstyle
\mathcal{S}^{(k)} = \set{x \in \mathcal{S} \st \text{$x \equiv y^{k} \mod{p}$
 for some $y$} }
\]
($k$-th powers $x$ modulo $p$) and
  \[
  \mathcal{S}_{d} = \set{ x \in \mathcal{S} \st \text{$x \equiv 0 \mod{d}$} }.
\]
  Then:

\begin{lemma} \label{mob1lemma}
    Let $\mathcal{S}$ be a set of integers and $e$ a divisor of $p-1$. Then
    \[
    \#\{x\in\mathcal{S} \st \gcd(x,p-1)=e\} = \sum_{k \mid \frac{p-1}{e}} \mu(k)\#\mathcal{S}_{e k},
    \]
    where $\mu(k)$ is the M\"obius function.
\end{lemma}

\begin{lemma}[Lemma~4 of~\cite{CZ}] \label{cz4lemma}
    Let $\mathcal{S}$ be a set of integers. Then
    \[
    \#\mathcal{S}^{\PR} = \sum_{k \mid p-1} \mu(k)\#\mathcal{S}^{(k)}.
    \]
\end{lemma}

\begin{lemma}[Lemma~5 of~\cite{CZ}] \label{cz5lemma}
    Let $p>3$ be a prime number, $\mathcal{P}=\mathcal{P}(a,r,N)$, and
    let $k$ and $d$ be integers between $1$ and $p-1$ such that $k$
    divides $p-1$.  Then
    \[
    \abs{\#\mathcal{P}^{(k)}_{d}- \frac{\#\mathcal{P}_{d}}{k}} \leq
    \sqrt{p}(1+\ln p).
    \]
\end{lemma}

It should be noted that~\cite{CZ} only proves Lemma~\ref{cz5lemma} for
$\gcd(r,d)=1$, but the proof goes through more generally.

Now the ``independence'' of $\RP$ and $\PR$:

  \begin{lemma}[Lemma~6 of~\cite{CZ}] \label{cz6lemma}
    Let
$\mathcal{P}=\mathcal{P}(a,r,N)$ with $\gcd(r,p-1)=1$. Then
\[
\abs{\#\mathcal{P}^{\RPPR} - N \left(\frac{\phi(p-1)}{p-1}\right)^{2}} \leq
\sigzero(p-1) + \sigzero(p-1)^{2}\sqrt{p}(1+\ln p).
\]
\end{lemma}

As the second author observed in~\cite{MoreeMR}, the factors of
$\sigzero(p-1)$ which occur here can in fact be improved to $\sum_{d
\mid p-1} \abs{\mu(d)} = 2^{\omega(p-1)}$ using the same proof; this
is also done in~\cite{Zhang}.  In addition, 
if $p-1 \mid N$ then the
first $\sigzero(p-1)$ term may be omitted.

In fact, several times in~\cite{Holden02} the following more
general heuristic was used:

\begin{heuristic} \label{indepheur}
    The order of $x$ modulo $p$ is independent of the
greatest common divisor of $x$ and $p-1$,  in the sense that for all $p$,
    \begin{multline*}
    \frac{1}{p-1}\#\set{x \in \set{1, \ldots, p-1} \st \gcd(x,p-1)=e, \quad
\ord_{p}(x)=\frac{p-1}{f} } \\
    \begin{aligned}
    &\approx \frac{1}{p-1}  \#\set{x \in \set{1, \ldots, p-1} \st \gcd(x,p-1)=e} \\
    &   \quad \times
    \frac{1}{p-1} \#\set{x \in \set{1, \ldots, p-1} \st \ord_{p}(x)=\frac{p-1}{f}}.
    \end{aligned}
    \end{multline*}
\end{heuristic}

To prove a rigorous form of this we need slightly less generality in
the sequence than in Lemma~\ref{cz6lemma}.  (The 
observations on Lemma~\ref{cz6lemma} likewise hold here.)

\begin{lemma} \label{ordgcdlemma}
    Let $e$ and $f$ be divisors of
$p-1$, and $N$ a multiple of $p-1$.  Let
$\mathcal{P}=\mathcal{P}(1,1,N)$ and
\[
\mathcal{P}' = \left\{ x \in \mathcal{P} \st \gcd(x,p-1)=e, \quad
\ord_{p}(x)=\frac{p-1}{f} \right\}.
\]
Then
\begin{multline*}
\begin{aligned}
    \abs{\#\mathcal{P}' - \frac{N}{(p-1)^2} \phi\left(\frac{p-1}{f}\right)
\phi\left(\frac{p-1}{e}\right)}
&\leq
 \sigzero \left(\frac{p-1}{f}\right)
\sigzero\left(\frac{p-1}{e}\right) \sqrt{p}(1+\ln p) \\
&\leq
 \sigzero\left(p-1\right)^{2}
 \sqrt{p}(1+\ln p).
\end{aligned}
\end{multline*}
\end{lemma}

  With the use of the more general version of Lemma~\ref{cz5lemma},
the proof of Lemma~\ref{ordgcdlemma} is essentially the same as that of
Lemma~\ref{cz6lemma}.

An equivalent way of thinking about Heuristic~\ref{indepheur} is to
fix a primitive root $b$ modulo $p$ and say that the discrete
logarithm $\log$ with base $b$ is a ``random map'' considered in terms
of divisibility; that is, that $\gcd(\log x, p-1)$ (which equals
$(p-1)/\ord_{p}(x)$) is distributed independently of $\gcd(x, p-1)$.
If we apply this discrete logarithm to~\eqref{fp}, we get a new equation:

\begin{equation} \label{fpprime}
    h \log g \equiv \log h \mod{p-1}.
\end{equation}

Looking at~\eqref{fpprime} with the ``random map'' idea in mind, we
see that $\gcd(g,p-1)$ seems to be independent of this equation.
This is the idea underlying the following heuristic:

\begin{heuristic} \label{grpheur}
    Among solutions to~\eqref{fp}, the greatest common divisor of
    $g$ and $p-1$ is independent of all other conditions on the order
    and greatest common divisor of $g$ and $h$, in the sense that for all $p$,
    \begin{multline*}
    \begin{aligned}
        &\frac{1}{p-1}\#\left\{g \st \text{\eqref{fp} holds,\ }
        \gcd(h,p-1)=e, \right. \\
        &\left. \phantom{\frac{1}{p-1}\#\{g \st}
        \ord_{p}(h)=\frac{p-1}{f}, \  \ord_{p}(g)=\frac{p-1}{d},
        \ \gcd(g,p-1)=n \right\} \\
    \end{aligned} \\
    \begin{aligned}
        &\approx \frac{1}{p-1}\#\set{g \st \text{\eqref{fp} holds,\ }
        \gcd(h,p-1)=e, \
        \ord_{p}(h)=\frac{p-1}{f}, \ \ord_{p}(g)=\frac{p-1}{d} } \\
        &   \quad \times
        \frac{1}{p-1} \#\set{g \st \text{\eqref{fp} holds,}
        \ \gcd(g, p-1)=n} .
    \end{aligned}
    \end{multline*}
\end{heuristic}

Heuristic~\ref{grpheur}, unlike Heuristics~\ref{rpprheur}
and~\ref{indepheur}, cannot yet be made rigorous.

\section{Conjectures for Fixed Points}

\label{fpconjsec}

The following conjectures and
theorems on fixed points were listed in~\cite{Holden02} and corrected
in the unpublished notes~\cite{Holden02a}.

    \begin{proposition} \label{prop1}
    $\Nfp_{g \ANY, h \RP}(p) = \phi(p-1).$
    \end{proposition}

    \begin{theorem}[Zhang, independently by Cobeli and Zaharescu]
    \label{thm1}
\begin{equation*}
\begin{split}
    \Nfp_{g \PR, h \RPPR}(p)& =  \Nfp_{g \PR, h \RP}(p) \\
    & = \Nfp_{g \PR, h \PR}(p) \\
    & = \Nfp_{g \ANY, h \RPPR}(p) \\
    & = \Nfp_{g \ANY, h \PR}(p) \\
    &\approx
\mbox{$\phi(p-1)^{2}/(p-1)$}.
\end{split}
\end{equation*}
\end{theorem}

    \begin{conjecture} \label{conj1}
    \mbox{}
    \begin{enumerate}
\item \label{conj1a} $\Nfp_{g \ANY, h \ANY}(p) \approx
p-1$.
\item \label{conj1b} $\Nfp_{g \PR, h \ANY}(p) \approx \phi(p-1)$.
\item \label{conj1d} $\Nfp_{g \RP, h \bullet}(p)
\approx \phi(p-1)/(p-1) \Nfp_{g \ANY, h \bullet}(p)$.
\item \label{conj1e} $\Nfp_{g \RPPR, h \bullet}(p)
\approx \phi(p-1)/(p-1) \Nfp_{g \PR, h \bullet}(p)$.
\end{enumerate}

\end{conjecture}

\begin{remark}
Note that Conjecture~1(c) of~\cite{Holden02} is incorrect.
In~\eqref{fp} if $h \PR$ then $g \PR$ also, so $\Nfp_{g \ANY, h
\PR}(p)$ is equal to $\Nfp_{g \PR, h \RPPR}(p)$ and not different as
was originally conjectured.
\end{remark}

Proposition~\ref{prop1} follows directly from the fact that $g =
h^{\overline{h}}$.  Theorem~\ref{thm1} also follows, with the
application of Lemma~\ref{cz6lemma}.  (That is,
Heuristic~\ref{rpprheur}.)  Conjecture~\ref{conj1}(\ref{conj1a}) is
essentially the same but we need to consider whether $h$ is an $e$-th
power, where $e=\gcd(h, p-1)$.  
Thus the conjecture uses Heuristic~\ref{indepheur}.  More
specifically, we see that~\eqref{fp} can be solved exactly when
$\gcd(h, p-1)=e$ and $h$ is an $e$-th power modulo $p$, and in fact
there are exactly $e$ such solutions.  Thus
\begin{equation} \label{NeTeqn}
\Nfp_{g \ANY, h \ANY}(p) = \sum_{e \mid p-1} e\ T(e,p).
\end{equation}
where
\[
T(e,p) =
\#\left\{ h \in
\mathcal{P}\left(1,1,{p-1}\right)^{(e)} \st \gcd(h,p-1)=e \right\}.
\]
According to Heuristic~\ref{indepheur}, we can model this sum using a 
set of independent random variables $X_{1}, \ldots, X_{p-1}$ such that
\[
X_{h} = \begin{cases}
    \gcd(h,p-1) & \text{with probability $\frac{1}{\gcd(h,p-1)}$}; \\
    0 & \text{otherwise.}
\end{cases}
\]
Then the heuristic suggests that $\Nfp_{g \ANY, h \ANY}(p)$ 
is approximately equal to the expected value of $X_{1} + \cdots + 
X_{p-1}$, which is clearly $p-1$.

Conjecture~\ref{conj1}(\ref{conj1b}) was justified in~\cite{Holden02}
using the argument that $g \PR$ should be independent of
 $\gcd(h,p-1)$
and $\ord_{p} h$.  This is somewhat dubious on the face of it, since
if~\eqref{fp} holds then the order of $g$ is certainly
constrained by both $\gcd(h,p-1)$ and $\ord_{p} h$.  The assumption is not
necessary, however.

Observe first that if~\eqref{fpprime} holds with $g
\PR$ then $\gcd(h,p-1) = \gcd(\log h, p-1)$.  Then we apply the
following elementary lemma:

\begin{lemma} \label{prsolslemma}
    Let $\gcd(a, q) = \gcd(b, q) = d$.  Then the number of solutions of
    \[
    a x \equiv b \mod{q}
    \]
    with $\gcd(x, q)=1$ is given by $\phi(q)/\phi(q/d)$.    
      In particular, there are always
    between $1$ and $d$ solutions.
\end{lemma}

Thus the number of solutions to~\eqref{fp} with $g \PR$ and $h \ANY$ is
\[
\sum_{d \mid p-1} \#\set{x \in \set{1, \ldots, p-1} \st \gcd(x,p-1)=d, \quad
\ord_{p}(x)=\frac{p-1}{d} } \frac{\phi(p-1)}{\phi\left((p-1)/d\right)}
\]
which by Heuristic~\ref{indepheur} is approximately equal to
\[
\sum_{d \mid p-1} \frac{1}{p-1}
\left(\phi\left(\frac{p-1}{d}\right)\right)^{2} \frac{\phi(p-1)}{\phi 
\left((p-1)/d\right)} = \phi(p-1).
\]
This argument justifies Conjecture~\ref{conj1}(\ref{conj1b}).

Conjectures~\ref{conj1}(\ref{conj1d}) and~\ref{conj1}(\ref{conj1e})
were justified in~\cite{Holden02} with Heuristic~\ref{grpheur};
in fact the conjectures are merely special cases of the heuristic.

In Section~\ref{fpthmsec}, we will try to approximate the error term
in Conjectures~\ref{conj1}(\ref{conj1a})
and~\ref{conj1}(\ref{conj1b}) using Lemma~\ref{ordgcdlemma}.  The
results, however, will not be entirely satisfactory.  With this in
mind, we will also use Heuristic~\ref{indepheur} to model the 
distribution of the values of $\Nfp_{g \ANY, h \ANY}(p)$.
Let $X_{1}, \ldots, X_{p-1}$ be as above.  Then we wish to find 
$\sigma^{2}$, the expected value of
\[
\left(\sum_{h=1}^{p-1} X_{h}-(p-1) \right)^{2}.
\]
Note that the expected value of $X_{h}X_{j}$ is $\gcd(h, p-1)$ if 
$h=j$ and $1$ otherwise.  Using this, an easy computation shows that
\begin{align*}
    \sigma^{2} &= \sum_{h=1}^{p-1} \gcd(h,p-1) - (p-1) 
     = \sum_{d \mid p-1} d \ \phi\left(\frac{p-1}{d}\right) - (p-1).
\end{align*}
In particular, $\sigma< p^{1/2+\epsilon}$ for every $\epsilon>0$.  
Thus we have the following:

\begin{conjecture} \label{varconj}
There are $o(x/\ln x)$ primes $p \leq x$ for which
\[
\abs{\Nfp_{g \ANY, h \ANY}(p) - (p-1)} > p^{1/2+\epsilon}
\]
for every $\epsilon>0$.
\end{conjecture}
Some progress toward proving this conjecture is described in 
Section~\ref{fpthmsec}.

Proposition~\ref{prop1}, Theorem~\ref{thm1}, and Conjecture~\ref{conj1}
are summarized in Table~\ref{fptalktable}, which appeared
in~\cite{Holden02a}.  The table also contains new data collected
since~\cite{Holden02}.

\begin{table}[!ht]
    \caption{Solutions to~\eqref{fp}}
    \label{fptalktable}
    $$\begin{array}{|l|l|l|l|l|}
\multicolumn{5}{l}{\text{(a) Predicted formulas for $\Nfp(p)$}}\\
    \hline
        g \setminus h & \ANY & \PR & \RP & \RPPR  \\
        \hline
    \ANY &\approx \scriptstyle(p-1) &
    \approx\frac{\phi(p-1)^{2}}{(p-1)} &
        = \scriptstyle \phi(p-1) &
    \approx\frac{\phi(p-1)^{2}}{(p-1)} \\

        \hline
    \PR & \approx \scriptstyle \phi(p-1) &
    \approx\frac{\phi(p-1)^{2}}{(p-1)} &
    \approx\frac{\phi(p-1)^{2}}{(p-1)} &
    \approx\frac{\phi(p-1)^{2}}{(p-1)} \\

    \hline
    \RP & \approx \scriptstyle\phi(p-1) &
    \approx\frac{\phi(p-1)^{3}}{(p-1)^2} &
    \approx\frac{\phi(p-1)^{2}}{(p-1)} &
    \approx\frac{\phi(p-1)^{3}}{(p-1)^2} \\

    \hline
    \RPPR & \approx\frac{\phi(p-1)^{2}}{(p-1)} &
    \approx\frac{\phi(p-1)^{3}}{(p-1)^2} &
    \approx\frac{\phi(p-1)^{3}}{(p-1)^2} &
    \approx\frac{\phi(p-1)^{3}}{(p-1)^2} \\
        \hline
\multicolumn{5}{l}{}\\
    \multicolumn{5}{l}{\text{(b) Predicted values for
    $\Nfp(100057)$}}\\
        \hline
        g \setminus h & \ANY & \PR & \RP & \RPPR  \\
        \hline
\ANY & 100056 & 9139.46 & 30240 & 9139.46 \\
    \hline
    \PR  &30240  &9139.46  &9139.46 & 9139.46 \\
    \hline
    \RP  &30240  &2762.23 & 9139.46 & 2762.23 \\
    \hline
    \RPPR  &9139.46&  2762.23 & 2762.23 & 2762.23  \\
    \hline
\multicolumn{5}{l}{}\\
    \multicolumn{5}{l}{\text{(c) Observed values for
    $\Nfp(100057)$}}\\

        \hline
        g \setminus h & \ANY & \PR & \RP & \RPPR  \\
        \hline
   \ANY&  98506& 9192 & 30240& 9192\\
   \hline
\PR&   29630& 9192&  9192&  9192\\
\hline
\RP&   29774& 2784&  9037&  2784\\
\hline
\RPPR& 9085&  2784&  2784&  2784\\
\hline
\end{array}$$
\end{table}

\section{Theorems on Fixed Points}
\label{fpthmsec}

The first rigorous result on this subject was Theorem~\ref{thm1}.
Both~\cite{Zhang} and~\cite{CZ} provided bounds on the error
involved; we will use notation closer to~\cite{CZ}.

\begin{theorem}[Theorem~1 of~\cite{CZ}] \label{cz1thm}
\[
\abs{\Nfp_{g \PR, h \RPPR}(p)  - \frac{\phi(p-1)^{2}}{p-1}} \leq
\sigzero(p-1)^{2}\sqrt{p}(1+\ln p).
\]
\end{theorem}

\begin{proof} Apply Lemma~\ref{cz6lemma} with
$\mathcal{P}=\mathcal{P}(1,1,p-1)$.  (The observations on
$\sigzero(p-1)$ apply.)
\end{proof}

We next turn our attention to $\Nfp_{g \ANY, h \ANY}(p)$.
Recall from Section~\ref{fpconjsec} that its value can be expressed
by~\eqref{NeTeqn}.  The quantity $T(e,p)$ which occurs there can be
straightforwardly evaluated using Lemmas~\ref{mob1lemma}
and~\ref{cz5lemma}.  We can also use the following characterization:

\begin{lemma}  \label{pieter4.5prime}
    Let $k \mid p-1$.
    Then
    \[
    T \left( \frac{p-1}{k}, p \right) = \#\set{ j \st 1 \leq j \leq 
    k, \quad (j,k)=1, \quad  (-j)^{k} \equiv k^{k} \mod{p} }.
    \]
\end{lemma}

\begin{proof}
    For each integer $h$ with 
    $\gcd(h,p-1) = (p-1)/k$,  $h=j (p-1)/k$ for some $1 \leq j 
    \leq k$ with $\gcd(j,k)=1$, such that, moreover, 
    \[
    j \ \frac{p-1}{k} \equiv x^{(p-1)/k} \mod{p}
    \]
    for some integer $x$.  It follows that
    \[
    \left( j \ \frac{p-1}{k} \right)^{k} \equiv 1 \mod{p}
    \]
    and hence
    \[
    (-j)^{k} \equiv k^{k} \mod{p}.
    \]
    (Note that $p \nmid k$.)  On observing that if 
    \[
    z^{k} \equiv 1 \mod{p},
    \]
    then 
    \[
    z \equiv x^{(p-1)/k} \mod{p}
    \]
    for some integer $x$, the proof of the reverse implication easily 
    follows.
\end{proof}

We now have the following results:

\begin{proposition} \label{anyanyprop}
Let $ e \mid p-1$.  Then
    \begin{enumerate}
	\item \label{Tpart} $\displaystyle \abs{T(e,p) - \frac{1}{e} \phi\left(\frac{p-1}{e} \right)}
	      \leq \sigzero \left(\frac{p-1}{e}\right) \sqrt{p}(1+\ln p).$
	\item \label{T1part} $\displaystyle T(1,p) = \phi(p-1).$
	\item \label{Tlargepart} If $k$ is a divisor of $p-1$ such that 
	$2k^{k} \leq p$, then 
	$\displaystyle T\left(\frac{p-1}{k},p\right) =0.$
	\item \label{Ttrivpart} $\displaystyle 0 \leq T(e,p) \leq 
	\phi\left(\frac{p-1}{e}\right).$
	\item \label{anyanypart} 
	$\displaystyle
	\abs{\Nfp_{g \ANY, h \ANY}(p)  - (p-3)} 
	\leq \sigzero(p-1)
	\left(\sigone(p-1)-\frac{3}{2}(p-1)\right) \sqrt{p}(1+\ln p).$
	\item \label{anyanypart2} For any $E$, $1 \leq E \leq p-1$,
	\[
	\abs{\Nfp_{g \ANY, h \ANY}(p)  - (p-1)} \leq E \ 
	\sigzero(p-1)^{2} \sqrt{p}(1+\ln p) + (p-1) 
	\sigzero_{\frac{p-1}{E}}(p-1),
	\]
	where
	\[
	\sigzero_{k}(n) = \#\set{d \mid (p-1) \st d < k}.  
	\]
    \end{enumerate}
\end{proposition}

\begin{proof}
The cardinality of $T(e,p)$ equals
    \begin{align*}
 &  \#\set{ h \in
\mathcal{P}\left(1,1,{p-1}\right)^{(e)} \st \gcd(h,p-1)=e }\\
&= \sum_{k \mid \frac{p-1}{e}} \mu(k) 
\#\mathcal{P}\left(1,1,{p-1}\right)^{(e)}_{ek} \\
\intertext{by Lemma~\ref{mob1lemma}}
& = \sum_{k \mid \frac{p-1}{e}} \mu(k) \left[\frac{1}{e}
\#\mathcal{P}\left(1,1,{p-1}\right)_{ek} + \eta_{e,k} \sqrt{p}(1+\ln p)
\right]\\
\intertext{for some $-1\leq \eta_{e,k} \leq 1$, by Lemma~\ref{cz5lemma}}
& = \sum_{k \mid \frac{p-1}{e}} \mu(k) \left[\frac{1}{e}
\frac{p-1}{e k} + \eta_{e,k} \sqrt{p}(1+\ln p)
\right] \\
&= \left[\left(\sum_{k \mid \frac{p-1}{e}}
\frac{\mu(k)}{k} \right) \frac{p-1}{e^2} +
\eta_e \sqrt{p}(1+\ln
p)\sigzero\left(\frac{p-1}{e}\right) \right] \\
\intertext{for some $-1 \leq \eta_e \leq 1$}
&= \frac{1}{e} \left[\phi\left(\frac{p-1}{e}\right) +
 \eta_e \sqrt{p}(1+\ln
p)\sigzero\left(\frac{p-1}{e}\right) \right] \\
\end{align*}
from whence part~\ref{Tpart} follows.  

Parts~\ref{T1part} and~\ref{Ttrivpart} are clear from the
definition.

Part~\ref{Tlargepart} follows from Lemma~\ref{pieter4.5prime}, since 
for such values of $k$ one has 
\[0 < k^{k} - (-j)^{k} < p\]
for any $j$ 
between $1$ and $k$, relatively prime to $k$.  (This was observed by 
an anonymous referee.)

Part~\ref{anyanypart} follows upon noting 
that
\begin{align*}
 \sum_{e \mid p-1} e \ T(e,p) 
 & = \sum_{e \mid p-1} \left[\phi\left(\frac{p-1}{e}\right) +
 e \ \eta_e \sqrt{p}(1+\ln p)\sigzero\left(\frac{p-1}{e}\right) \right] \\
&= (p-1) +
  \eta \sqrt{p}(1+\ln p)
\sigzero\left(p-1\right)\sigone(p-1)  \\
\end{align*}
for some $-1 \leq \eta \leq 1$ and then applying part~\ref{Tlargepart}.

Part~\ref{anyanypart2} is similar; observe that
\begin{multline*}
    \sum_{e \mid p-1} e \ T(e,p) \\
\begin{aligned}
 & = \sum_{\substack{e \mid p-1 \\ e \leq E}}
 \left[\phi\left(\frac{p-1}{e}\right) +
 e \ \eta_e \sqrt{p}(1+\ln p)\sigzero\left(\frac{p-1}{e}\right) \right] 
 + \sum_{\substack{e \mid p-1 \\ e > E}} e \ T(e,p) \\
 &= \sum_{\substack{e \mid p-1 \\ e \leq E}}
 \left[\phi\left(\frac{p-1}{e}\right) +
 e \ \eta_e \sqrt{p}(1+\ln p)\sigzero\left(\frac{p-1}{e}\right) \right] 
 + \eta' \sum_{\substack{e \mid p-1 \\ e > E}} e \ 
 \phi\left(\frac{p-1}{e}\right) \\
 & = (p-1) + \sum_{\substack{e \mid p-1 \\ e \leq E}}
 e \ \eta_e \sqrt{p}(1+\ln p)\sigzero\left(\frac{p-1}{e}\right) 
 + \eta' \sum_{\substack{e \mid p-1 \\ e > E}} (e-1) 
 \phi\left(\frac{p-1}{e}\right) \\
 &=  (p-1) + E\ \eta\ \sigzero(p-1)^{2}
   \sqrt{p}(1+\ln p) +  \eta' \sum_{\substack{e \mid p-1 \\ e >
 E}} (p-1) \\
 &= (p-1) + E\ \eta\ \sigzero(p-1)^{2} \sqrt{p}(1+\ln p) + \eta' (p-1) 
 \sigzero_{\frac{p-1}{E}}(p-1),
\end{aligned}
\end{multline*}
where $-1 \leq \eta \leq 1$, $-1 \leq \eta_{e} \leq 1$, $-1 \leq
\eta' \leq 1$.

\end{proof}

Unfortunately for part~\ref{anyanypart}, $\sigone(p-1)-3(p-1)/2 = O(p \ln \ln p)$ in the worst
case, although if $p$ is a Sophie Germain prime $\sigone(p-1)-3(p-1)/2
=3$, and the ``average case'', averaging over a range of $p$, is
$\sigone(p-1)-3(p-1)/2 \approx 0.70386 (p-1)$.  (See later in this
section for more on Sophie Germain primes, and
Sections~\ref{aveconjsec} and~\ref{averessec} for further details of
the ``average case''.)  Thus the ``error'' term for $\Nfp_{g \ANY, h
\ANY}(p)$ will be larger than the main term for infinitely many $p$.
In fact, this estimate is even weaker than the rather trivial bound
\[
\phi(p-1) \leq \Nfp_{g \ANY, h \ANY}(p) 
\leq \sum_{e \mid p-1} e \ \phi\left(\frac{p-1}{e}\right) \leq (p-1) \sigzero(p-1)
\]
obtained from parts~\ref{T1part} and~\ref{Ttrivpart} of the
proposition.  (On the basis of an heuristic argument we conjecture that
the average order of
\[
\sum_{e \mid p-1} e \ \phi\left(\frac{p-1}{e}\right)
\]
is $c_{1} p \ln p$ with $c_{1}$ a positive constant.)
A little thought reveals the problem: since $\#\{ h \in
\mathcal{P}(1,1,p-1)^{(e)} \st \gcd(h,p-1)=e \}$ is multiplied by
each divisor $e$ of $p-1$; an error of even $1$ in calculating the
number of elements in the set for a large value of $e$ will result in
an error of $O(p-1)$.

Part~\ref{anyanypart2} gives us something of an improvement; but
it does not solve the problem in general.  In order to make the term
$E\ \sigzero(p-1)^{2} \sqrt{p}(1+\ln p)$ be even $O(p-1)$, we must 
pick $E < \sqrt{p-1}$, which makes $\sigzero_{\frac{p-1}{E}}(p-1) \leq  
\sigzero(p-1)/2$ by elementary counting of divisors.  Thus the 
``error'' term will still be of larger order than the main term.

On the other hand, the line of argument from
part~\ref{anyanypart} works if we restrict to primes $p$ for which
\[
E(p) = \max \set{ e \st e \mid p-1, \quad T(e,p) > 0 }
\]
is not too large.  (Thus, the error in $T(e,p)$ will not be multiplied 
by too large an $e$.)

\begin{proposition} \label{pieter4.3}
    Suppose that $1/4 \leq \beta \leq 1$, $E(p) \leq p^{\beta}$, 
    and $\delta > 0$ then 
    \[
    \Nfp_{g \ANY, h \ANY}(p)  = (p-1) + O_{\delta}\left( 
    p^{1/2+\beta+\delta}\right).
    \]
    More specifically,
    \[
       \abs{\Nfp_{g \ANY, h \ANY}(p) - (p-1)} \leq
       p^{1/2+\beta} \sigzero(p-1)^{2} (2 + \ln p).       
    \]
\end{proposition}

\begin{proof}
    By the assumption on $E(p)$, \eqref{NeTeqn}, and 
    Proposition~\ref{anyanyprop}(\ref{Tpart}), we have:
    
    \begin{align*}
	\Nfp_{g \ANY, h \ANY}(p) &= \sum_{\substack{e \mid p-1 \\ 
	e \leq p^{\beta}}} e \ T(e,p) \\ 
	& = \sum_{\substack{e \mid p-1 \\ e \leq p^{\beta}}} 
	\phi\left( \frac{p-1}{e} \right) +
	\eta_{1} p^{\beta} \sigzero(p-1)^{2} \sqrt{p} (1 + \ln p)  \\
	\intertext{for some $-1 \leq \eta_{1} \leq 1$}
	& = p-1 - \left( \sum_{\substack{e \mid p-1 \\ e > p^{\beta}}}
	\phi\left( \frac{p-1}{e} \right) \right) + \eta_{1}
	p^{1/2+\beta} \sigzero(p-1)^{2}  (1 + \ln p)  \\
	& = p-1 + \eta_{2} \sigzero(p-1)  p^{1-\beta}+ \eta_{1}
	p^{1/2+\beta} \sigzero(p-1)^{2}  (1 + \ln p)  \\
	\intertext{for some $-1 \leq \eta_{2} \leq 0$}
	& = p-1 + \eta_{3} 
	p^{1/2+\beta} \sigzero(p-1)^{2}  (2 + \ln p)  \\
	\intertext{for some $-1 \leq \eta_{3} \leq 1$}
	& = p-1 + 
	O\left(p^{1/2+\beta+\delta} \right),
    \end{align*}
    where we used the facts that $\sigzero(n) = O_{\delta} 
    \left(n^{\delta}\right)$ for every $\delta > 0$ and $\phi(n) \leq n$.
\end{proof}

\begin{remark} One reason to consider the more specific version of this 
proposition is to aid in computer searches such as the one described 
in~\cite{CampbellThesis}.
\end{remark}

Proposition~\ref{pieter4.3} is, of course, only useful if there exist 
sufficiently many primes satisfying $E(p) \leq p^{\beta}$ for some 
appropriate $\beta$.  For instance, $\beta$ needs to be less than 
$1/2$ before the error term is less than the main term:

\begin{corollary}
    Suppose $E(p) \leq p^{1/2-\delta}$ 
    and $\delta > 0$. Then 
    \[
    \Nfp_{g \ANY, h \ANY}(p)  = (p-1) + o(p).
    \]
\end{corollary}

In fact, we will prove that there are $\gg x/\ln x$ primes $p \leq x$
for which $E(p) \leq p^{0.3313}$ and thus that there are $\gg x/\ln x$
primes $p \leq x$ such that
\[
\Nfp_{g \ANY, h \ANY}(p)  = (p-1) + O\left( p^{5/6} \right).
\]
The proof of this starts with the 
following application of Lemma~\ref{pieter4.5prime}:

\begin{proposition} \label{pieter4.6}
    Let $\delta > 0$, $\alpha \geq 2/3$.
    Except for 
    $
    O\left( x^{3-3\alpha}/{\ln^{3\alpha-1+3\delta} x} \right)
    $
    primes $p \leq x$ we have 
    \[
    E(p) < p^{\alpha} \ln^{\alpha+\delta} p.
    \]
    In particular, letting $\alpha=2/3$, except for
    $
    O\left({x}/{\ln^{1+3\delta} x} \right)
    $
    primes $p \leq x$ we have 
    \[
    E(p) < p^{2/3} \ln^{2/3+\delta} p.
    \]
\end{proposition}

\begin{proof}
    Let $f_{\delta}(x) = x^{1-\alpha}/\ln^{\alpha+\delta} x$.  If $p \leq 
    x$ is a prime not dividing
    \[
    P = \prod_{1\leq k \leq f_{\delta}(x)} \prod_{\substack{j=1 \\ 
    (j,k)=1}}^{k} \left( (-j)^{k} - k^{k} \right),
    \]
    then, by Lemma~\ref{pieter4.5prime},
    \[
    T\left( \frac{p-1}{k_{1}}, p \right) > 0
    \]
    for some $k_{1}$ implies
    \[
    k_{1} > f_{\delta}(x) > \frac{p-1}{p^{\alpha}\ln^{\alpha+\delta} p}
    \]
    and hence
    \[
    E(p) < p^{\alpha} \ln^{\alpha+\delta} p.
    \]    
    The non-zero integer $P$ has 
    \[
    O\left(\sum_{k \leq f_{\delta}(x)} k^{2} \ln k \right)
    = O\left(f_{\delta}(x)^{3} \ln f_{\delta}(x) \right) = 
    O \left ( \frac{x^{3-3\alpha}}{\ln^{3\alpha-1+3\delta} x} \right)
    \]
    distinct prime divisors.  These are the possible exceptions to 
    the inequality
        \[
    E(p) < p^{\alpha} \ln^{\alpha+\delta} p.
    \]
\end{proof}

We can now prove:

\begin{proposition} \label{pieter4.4}
 There are $\gg x/\ln x$ primes $p \leq x$ 
for which $E(p) \leq p^{0.3313}$.
\end{proposition}

\begin{proof} It is a deep result of Fouvry (see, e.g., 
\cite{Fouvry}), that $\gg x/\ln x$ primes $p \leq x$ are such that 
$p-1$ has a prime factor larger than $p^{0.6687}$.  In combination 
with Proposition~\ref{pieter4.6} it follows that there are $\gg x/\ln 
x$ primes $p \leq x$ for which $E(p) < p^{0.668}$ and $p-1$ has a prime 
factor larger than $p^{0.6687}$.  Since $E(p)$ is a divisor of $p-1$ 
it must divide the factors of $p-1$ besides the largest, and thus
$E(p) < p^{0.3313}$ for any such primes.
\end{proof}

Letting $\beta=0.3313$ and $\delta=0.002$ in
Proposition~\ref{pieter4.3} and invoking Proposition~\ref{pieter4.4},
we now have:

\begin{theorem} \label{pieterA}
    There are $\gg x/\ln x$ primes $p \leq x$ such that
        \[
    \Nfp_{g \ANY, h \ANY}(p)  = (p-1) + O\left( 
    p^{5/6}\right).
    \]
    More specifically, there are $\gg x/\ln x$ primes $p \leq x$ such that
    \[
       \abs{\Nfp_{g \ANY, h \ANY}(p) - (p-1)} \leq
       p^{0.8313} \sigzero(p-1)^{2} (2 + \ln p).       
    \]
\end{theorem}

\begin{remark} 
    If one can establish that in Fouvry's assertion, 0.6687
can be replaced by some larger $\theta$ (up to $\theta = 3/4$), then in
Theorem~\ref{pieterA} the exponents $5/6$ and $0.8313$ can be replaced
by $3/2-\theta+\delta$ and $3/2-\theta$ for any $\delta>0$.
\end{remark}

The most well-known primes $p$ with $p-1$ having a large prime factor 
are the Sophie Germain primes.  These are the primes $p$ such that 
$p-1=2q$ with $q$ a prime.  For these primes it is easily shown 
(using Proposition~\ref{anyanyprop}(\ref{Tlargepart}) with $k=1$ and 
$k=2$) that
\[
\Nfp_{g \ANY, h \ANY}(p) = T(1,p) + 2T(2,p).
\]
Proceeding as in the proof of 
Proposition~\ref{anyanyprop}(\ref{anyanypart}), the following
result is then obtained:

\begin{proposition} If $p$ is a Sophie Germain prime, then
   \[
   \abs{\Nfp_{g \ANY, h \ANY}(p) - (p-3)} \leq 2\sqrt{p} (1+\ln 
   p)
   \]
\end{proposition}

By sieving methods it can be shown that there are $\ll x/\log^{2} x$ 
Sophie Germain primes $p \leq x$.  On the other hand, it is not known 
whether or not there are infinitely many Sophie Germain primes.

In fact, we can state a similar result for primes $p$ of the form 
$p-1=mq$ as long as $q$ is prime and $m$ is sufficiently small.  (This 
was observed by an anonymous referee.)  Let $W$ be the Lambert $W$ 
function, which has the property that $W(x)e^{W(x)} = x$ for any 
$x$.  Then as long as $m \leq \ln(p/2) / W(\ln(p/2))$, any divisor $k$ 
of $m$ will have the property that $2k^{k} \leq p$.  Thus
Proposition~\ref{anyanyprop}(\ref{Tlargepart}) gives us
\[
\Nfp_{g \ANY, h \ANY}(p) = \sum_{e \mid m} e \ T(e,p)
\]
and thus:

\begin{proposition} If $p$ is a prime as described above, then
   \[
   \abs{\Nfp_{g \ANY, h \ANY}(p) - (p-1-m)} \leq 2
   \sigzero(m)\sigone(m) \sqrt{p} (1+\ln p).
   \]
\end{proposition}

(The factor $2 \sigzero(m)\sigone(m)$ can sometimes be improved, as 
was the case for Sophie Germain primes.)

It is also worth asking how large $E(p)$ can be with respect to $p$.  
We put forward the following conjecture:
\begin{conjecture} Let $\alpha<1$.  There exist infinitely many 
primes $p$ with $E(p)>p^\alpha$.
\end{conjecture}

The idea is that amongst the numbers of the form 
\[
k^{k}-(-j)^{k}, \quad 1 \leq j \leq k, \quad \gcd(j,k)=1,
\]
there will be many that are close to being a prime and that if $q$ is 
a large prime divisor of such a number, then $E(q)$ will be large.  
Taking $k=29$ and $j=5$ we infer, for example, that the prime
\[
q=\frac{29^{29}+5^{29}}{34}
\]
satisfies $E(q) > q^{0.964}$.  If $k$ is odd and 
\[
q=k^{k}-(-j)^{k}
\] 
is a prime for some $1 \leq j \leq k$, then $E(q) > (q-1)^{1-1/k}$.

Turning back to the general case, the situation
where $g$ is $\PR$ and $h$ is $\ANY$ follows the argument
explained in the justification of
Conjecture~\ref{conj1}(\ref{conj1b}), and uses
Lemma~\ref{ordgcdlemma} to estimate the error term.  It is very
similar to the previous case, and unfortunately has the same problem 
in the general case:

\begin{proposition} \label{pranyprop} \mbox{}
    \begin{enumerate} 
	\item 
	\[
	\abs{\Nfp_{g \PR, h \ANY}(p)  - \phi(p-1)-2} 
	\leq
	\sigzero(p-1)^{2}
	\left(\sigone(p-1)-\frac{3}{2}(p-1)\right) \sqrt{p}(1+\ln p).
	\]
	\item 
	 For any $E$, $1 \leq E \leq p-1$,
		\[
		\abs{\Nfp_{g \PR, h \ANY}(p)  - \phi(p-1)} \leq E \ 
		\sigzero(p-1)^{2} \sqrt{p}(1+\ln p) + \phi(p-1) 
		\sigzero_{\frac{p-1}{E}}(p-1).  
		\]
	    \end{enumerate}
\end{proposition}

We can proceed in the same fashion as Theorem~\ref{pieterA}, however, 
to prove:

\begin{theorem} 
    There are $\gg x/\ln x$ primes $p \leq x$ such that
        \[
	\Nfp_{g \PR, h \ANY}(p)  = \phi(p-1) + O\left( 
    p^{5/6}\right).
    \]
    More specifically, there are $\gg x/\ln x$ primes $p \leq x$ such that
    \[
       \abs{\Nfp_{g \PR, h \ANY}(p) - \phi(p-1)} \leq
       p^{0.8313} \sigzero(p-1)^{3} (2 + \ln p).       
    \]
\end{theorem}

Finally, we should mention that the second author (in~\cite{MoreeMR})
pointed out that we could also estimate the number $\Mfp_{g
\PR, h \ANY}(p)$ of values $h$ such that there exists \emph{some} $g$
satisfying~\eqref{fp}, with $g \PR$ and $h \ANY$.  From
\[
 \Mfp_{g \ANY, h \ANY}(p) = \sum_{e \mid p-1} \#\set{ h \st
 \ord_{p} h = \frac{p-1}{e}, \ \gcd(h,p-1)=e},
\]
it was shown:

\begin{theorem} \label{thm5}
    \[
\abs{\Mfp_{g \PR, h \ANY}(p)  - \frac{1}{p-1} \sum_{e \mid
p-1} \phi \left( \frac{p-1}{e} \right)^{2}} \leq
\sigzero(p-1)^{3}\sqrt{p}(1+\ln p).
\]
\end{theorem}

Similarly, we can estimate
\[
\Mfp_{g \ANY, h \ANY}(p)
 = \sum_{e \mid p-1} \#\set{ h \in
\mathcal{P}(1,1,p-1)^{(e)} \st \gcd(h,p-1)=e },
\]
giving:

\begin{theorem} \label{thm6}
    \[
\abs{\Mfp_{g \ANY, h \ANY}(p)  - \sum_{e \mid
p-1}\frac{1}{e}  \phi \left( \frac{p-1}{e} \right)} \leq
\sigzero(p-1)^{2}\sqrt{p}(1+\ln p).
\]
\end{theorem}

Since we are no longer counting multiple solutions for each value of
$h$ the problem with the error terms discussed above disappears; the
error terms are $O(p^{1/2+\epsilon})$ while the main terms look on
average like a constant times $p$.

(For completeness, we should note that if $h$ is $\RP$ and/or $\PR$,
then
\[
\Mfp_{g \bullet, h \bullet}(p) = \Nfp_{\eqref{fp}, g \bullet, h
\bullet}(p).
\]
Heuristic~\ref{grpheur} would also predict that
\[
\Mfp_{g \RP, h \bullet}(p)
\approx \phi(p-1)/(p-1) \Mfp_{g \ANY, h \bullet}(p)
\]
and
\[
\Mfp_{g \RPPR, h \bullet}(p)
\approx \phi(p-1)/(p-1) \Mfp_{g \PR, h \bullet}(p).)
\]

\section{Equivalence of the Equations for Two-cycles}

\label{equivsec}

As observed in~\cite{Holden02}, conditions on~\eqref{tc} can sometimes
be translated into conditions on~\eqref{ha} in a relatively
straightforward manner.  Table~\ref{relationstable}, reproduced
from~\cite{Holden02a}, summarizes these straightforward relationships.

\begin{table}
    \caption{Relationship between solutions to~\eqref{tc}
and solutions to~\eqref{ha}}
    \label{relationstable}
$$\begin{array}{|l|l|l|l|l|}
    \hline
        a \setminus h & \ANY & \PR & \RP & \RPPR  \\
        \hline
    \ANY & & & g \ANY & g \PR \\
     & & &  h \RP &  h \RPPR \\
        \hline
    \PR & & & g \PR & g \PR\\
        & & &  h \RP &  h \RPPR \\
    \hline
    \RP & h \ANY & g \PR & h \RP & g \PR \\
    &  g \ORD h &  h \PR &  g \ORD h & h \RPPR \\

    \hline
    \RPPR & g \PR & g \PR & g \PR & g \PR\\
    &  h \RPPR &  h \RPPR &  h \RPPR &  h \RPPR\\
        \hline

    \end{array}$$
\end{table}

We can go slightly further, however.  Taking the logarithm
of the two equations of~\eqref{tc} with respect to the same primitive root
$b$ gives us new equations:

\begin{equation} \label{tcprime}
    \begin{split}
    h \log g &\equiv \log a \mod{p-1}; \\
    a \log g &\equiv \log h \mod{p-1}.
    \end{split}
\end{equation}

Let $d=\gcd(h,a,p-1)$, and let $u_{0}$ and $v_{0}$ be such that
\[
u_{0} h + v_{0} a \equiv d \mod{p-1}.
\]
By using the Smith Normal Form, we can show that~\eqref{tcprime} is
equivalent to the equations:

\begin{equation}
    \begin{split}
    0 &\equiv \frac{h}{d} \log h - \frac{a}{d} \log a \mod{p-1}; \\
    d \log g &\equiv v_{0} \log h + u_{0} \log a \mod{p-1},
    \end{split}
\end{equation}
or:

\begin{equation} \label{tcsmith}
    \begin{split}
    h^{h/d} &\equiv a^{a/d} \mod{p}; \\
    g^{d}&\equiv h^{v_{0}} a^{u_{0}} \mod{p}.
    \end{split}
\end{equation}
In the case where $d=\gcd(h,a,p-1)=1$ then this becomes just

\begin{equation}
    \begin{split}
    h^{h} &\equiv a^{a} \mod{p}; \\
    g & \equiv h^{v_{0}} a^{u_{0}} \mod{p}.
    \end{split}
\end{equation}

Thus:

\begin{proposition} \label{equivprop}
    If $\gcd(h,a,p-1)=1$, then there is a one-to-one correspondence
    between triples $(g,h,a)$ which satisfy~\eqref{tc} and pairs
    $(h,a)$ which satisfy~\eqref{ha}, and the value of $g$ is unique
    given $h$ and $a$.  In particular, this is true if $h$ is $\RP$ or
    $a$ is $\RP$.
\end{proposition}

It was observed in~\cite{Holden02} that when neither $h$ nor $a$ is
$\RP$ the relationship between~\eqref{tc} and~\eqref{ha} is less
clear.  It was claimed there that given a pair $(h,a)$ which is a
solution to~\eqref{ha} we expect on average
$\gcd(a,p-1)\gcd(h,p-1)/\gcd(ha, p-1)^{2}$ pairs $(g,h)$ which are
solutions to~\eqref{tc}.

It is clear from~\eqref{tcsmith}, however, that when $d =
\gcd(h,a,p-1) \neq 1$ this is not the correct way to think about
things.  The proper equation to look at in this case is
not~\eqref{ha}, but
\begin{equation} \label{hasmith}
    h^{h/d} \equiv a^{a/d} \mod{p}.
\end{equation}
We will use $\Nha'$ to denote the number of solutions to~\eqref{hasmith}.

Now~\eqref{tcsmith} shows that a nontrivial solution to~\eqref{hasmith}
produces $d$ pairs $(g,h)$ which are nontrivial solutions to~\eqref{tc}
if $h^{v_{0}} a^{u_{0}}$ is a $d$-th power modulo $p$, and otherwise no
solutions.  (As in~\cite{Holden02}, we consider the ``trivial''
solutions to~\eqref{tc} to be the ones that are also solutions
to~\eqref{fp}.)  Thus the following heuristic implies that every
nontrivial solution to~\eqref{hasmith} produces on average one pair
$(g,h)$ which is a nontrivial solution to~\eqref{tc}.

\begin{heuristic} \label{uvheur}
For any pair $(h,a)$, let $d=\gcd(h,a,p-1)$, and let $u_{0}$ and
$v_{0}$ be such that
\[
u_{0} h + v_{0} a \equiv d \mod{p-1}.
\]
Then $(h,a) \mapsto h^{v_{0}}a^{u_{0}}$ is a random map even when
restricted to $\gcd(h,a,p-1)=d$, in the sense that
\begin{equation*}
\frac{\#\set{(h,a) \st h^{v_{0}}a^{u_{0}} \equiv y
\mod{p}, \ \gcd(h,a,p-1)=d} }%
{\#\set{(h,a) \st \gcd(h,a,p-1)=d} }
\approx
\frac{1}{\#\set{y \in \set{1, \ldots, p-1}}}.
\end{equation*}
\end{heuristic}

On the other hand, there is a solution to~\eqref{tcsmith} with $g \PR$
if and only if $h^{v_{0}} a^{u_{0}}$ is exactly a $d$-th power modulo
$p$; that is, $\ord_{p} (h^{v_{0}} a^{u_{0}}) = (p-1)/d$.  Then
Lemma~\ref{prsolslemma} says that the number of such solutions is
$\phi(p-1)/\phi((p-1)/d)$.
Thus Heuristic~\ref{uvheur} implies that every
solution to~\eqref{hasmith} produces on average
$\phi(p-1)/(p-1)$
pairs $(g,h)$ which
are solutions to~\eqref{tc} with $g \PR$.  These relationships 
between conditions on~\eqref{tc} and conditions on~\eqref{hasmith} 
are summarized in Table~\ref{relationstable2}, where 
$\mathbb{E}(\Ntc /\Nha')$ is the expected 
number of solutions to~\eqref{tc} given a solution to~\eqref{hasmith}.

    \begin{table}
    \caption{Relationship between solutions to~\eqref{tc}
and solutions to~\eqref{hasmith}}
    \label{relationstable2}

    $$\begin{array}{|c|c|c|c|c|}
    \hline
        g \setminus h & \ANY & \PR & \RP & \RPPR  \\
        \hline
    \ANY & h \ANY, a \ANY & h \PR & h \RP & h \RPPR \\
     & \mathbb{E}(\Ntc /\Nha')
     \approx 1 & a \RP & a \ANY &  a \RPPR \\
        \hline
    \PR & h \ANY, a \ANY & h \PR & h \RP & h \RPPR \\
     & \mathbb{E}(\Ntc /\Nha')
     \approx \frac{\phi\left({p-1}\right)}{p-1} & a \RP &
     a \PR & a \RPPR \\
        \hline

    \end{array}$$
    \end{table}

\section{Heuristics and Conjectures for Two-Cycles}

\label{tcconjsec}

We mentioned in Section~\ref{indepsec} that we could view $ x \mapsto \log
x$ as a ``random map'' in some sense.
We will also suppose that the map $x \mapsto x^{x} \mod{p}$
is ``random'', in a slightly different sense.

\begin{heuristic}
    The map $x \mapsto x^{x} \mod{p}$ is a random map given the
    obvious restrictions on order, in the sense that
    for all $p$, given $y \in \set{1, \ldots, p-1}$,
    then
    \begin{multline*}
    \#\set{x \in \set{1, \ldots, p-1} \st x^{x} \equiv y \mod{p} } \\
    \approx \frac{ \#\set{z \in \set{1, \ldots, p-1} \st
    (\ord_{p} z)/\gcd(z, \ord_{p} z) = \ord_{p} y } }%
    {\#\set{w \in \set{1, \ldots, p-1} \st \ord_{p} w = \ord_{p} y}}.
    \end{multline*}
\end{heuristic}
(The fraction on the right-hand side was referred to
in~\cite{Holden02} as ${\# S_{m}}/{\# T_{m}}$, where
$m=\ord_{p} y$.  The arguments there used this heuristic
implicitly.)

In fact, we would like a slightly stronger version of this:
\begin{heuristic} \label{xtoxxheur}
    The map $x \mapsto x^{x} \mod{p}$ is a random map even when restricted
    to a specific order and greatest common divisor,
     in the sense that
    for all $p$, given $y \in \set{1, \ldots, p-1}$ such that $\ord_p y =
    f/\gcd(e,f)$,
    then
    \begin{multline*}
    \#\set{x \in \set{1, \ldots, p-1} \st x^{x} \equiv y \mod{p},
    \ \gcd(x,p-1)=e, \ \ord_p x = f } \\
    \approx \frac{ \#\set{z \in \set{1, \ldots, p-1} \st
    \gcd(z,p-1)=e, \ \ord_p z = f } }%
    {\#\set{w \in \set{1, \ldots, p-1} \st \ord_{p} w = \ord_{p} y}}.
    \end{multline*}

\end{heuristic}

Heuristic~\ref{xtoxxheur}, like Heuristic~\ref{grpheur}, cannot yet
be made rigorous.

Using Proposition~\ref{equivprop} and Heuristic~\ref{xtoxxheur}, we
have the following conjectures from~\cite{Holden02}, as corrected
in~\cite{Holden02a}.

\begin{conjecture}\mbox{}
\label{conj23456}
\begin{enumerate}
\item     \label{conj2}
\(
    \Ntc_{g \ANY, h \RP}(p) = \Nha_{h \RP, a \ANY}(p) \approx 2
    \phi(p-1).
\)

\item \label{conj3}
\(
 \Ntc_{h \RP, g \ORD h}(p) = \Nha_{h \RP, a \RP}(p) \approx \phi(p-1)
 + \mbox{$\phi(p-1)^{2}/(p-1)$}.
\)

\item \label{conj4}
\(
   \Ntc_{g \PR, h \RP}(p) = \Nha_{h \RP, a \PR}(p) \approx
   2\phi(p-1)^{2}/(p-1).
\)

\item \label{conj5}
         \begin{equation*}
    \begin{split}
    \Ntc_{g \PR, h \RPPR}(p) &= \Ntc_{g \ANY, h \RPPR}(p) \\
&= \Nha_{h \RPPR, a \bullet}(p) 
= \Nha_{h \bullet, a \RPPR}(p) 
\\
&\approx
\phi(p-1)^{2}/(p-1) + \phi(p-1)^{3}/(p-1)^{2}.
\end{split}
\end{equation*}

\item \label{conj6a}
\(
\Ntc_{h \ANY, g \ORD h}(p) = \Nha_{h \ANY, a \RP}(p) \approx 2
\phi(p-1).
\)

\item \label{conj6b}
\(
\Ntc_{g \PR, h \PR}(p) = \Ntc_{g \ANY, h \PR}(p) = \Nha_{h \PR, a
\RP}(p) \approx 2\phi(p-1)^{2}/(p-1).
\)
\end{enumerate}
\end{conjecture}

In Conjectures~\ref{conj23456}(\ref{conj2})
and~\ref{conj23456}(\ref{conj6a}) it should be noted that the observed
values in question must be exactly (not just approximately) equal, by
the symmetry of~\eqref{ha}.  The same applies in
Conjectures~\ref{conj23456}(\ref{conj4})
and~\ref{conj23456}(\ref{conj6b}).

We also made in~\cite{Holden02} the following conjectures about
solutions to~\eqref{ha}.

\begin{conjecture} \label{conj7} \mbox{}
    \begin{enumerate}
\item \label{conj7a} $\displaystyle \Nha_{h \ANY, a \ANY}(p) \approx
(p-1) + \sum_{m \mid p-1} \frac{\phi(m)}{m^{2}} \left( \sum_{d \mid
(p-1)/m} \frac{\phi(dm)}{d}\right)^{2}$.

\item \label{conj7b} If $p-1$ is squarefree then $\displaystyle \Nha_{h
\ANY, a \ANY}(p) \approx (p-1) + \prod_{q \mid p-1} \left( q + 1 -
\frac{1}{q} \right)$, where the product is taken over primes $q$
dividing $p-1$.  

\item \label{conj7c} In general,
    \begin{multline*} \label{geneq}
    \Nha_{h \ANY, a \ANY}(p)  \\
\shoveleft{\approx (p-1) +  \prod_{q^\alpha \Vert p-1} \left(
\left[\left(1-\frac{1}{q}\right)\alpha + 1\right]^2 \right.} \\
+ \left(1-\frac{1}{q}\right)^3
\left[ (\alpha+1)^2 \frac{q^{\alpha+1}-q}{q-1}
- 2 (\alpha+1) \frac{\alpha q^{\alpha+2} -
(\alpha+1)q^{\alpha+1}+q}{(q-1)^2} \right.\\
\left.\left.
+ \frac{\alpha^2 q^{\alpha+3} - (2\alpha^2 + 2\alpha-1)q^{\alpha+2}
        + (\alpha^2+2\alpha+1)q^{\alpha+1}-q^2-q}{(q-1)^3}
\right]
\right)
\end{multline*}
where the product is taken over primes $q$ dividing $p-1$ and $\alpha$ is
the exact power of $q$ dividing $p-1$.

\item \label{conj7d} $\Nha_{h \PR, a \ANY}(p) \approx 2 \phi(p-1)$.

\item \label{conj7e}
$\Nha_{h \ANY, a \PR}(p)\approx 2 \phi(p-1)$.

\item \label{conj7f}
 $\Nha_{h \PR, a \PR}(p) \approx \phi(p-1) + \phi(p-1)^{2}/(p-1)$.

\end{enumerate}
\end{conjecture}

(The formulas in Conjecture~\ref{conj7}(\ref{conj7a}) and
Conjecture~\ref{conj7}(\ref{conj7c})
appear in~\cite{Holden02} with typos.  They appear correctly here and
in~\cite{Holden02a}.)

These conjectures rely on Heuristics~\ref{indepheur}
and~\ref{xtoxxheur} and a standard birthday paradox argument.  Thanks
to Lemma~\ref{ordgcdlemma} we are now closer to making them into
rigorous theorems.  All of the conjectures on~\eqref{ha} are summarized
in Table~\ref{hatalktable}, which appeared in~\cite{Holden02a}.  The
table also contains new data collected since~\cite{Holden02}.  As
in~\cite{Holden02}, we distinguish between the ``trivial'' solutions
to~\eqref{ha}, where $h=a$, and the ``nontrivial'' solutions.

\begin{table}[!ht]
    \caption{Solutions to~\eqref{ha}}
    \label{hatalktable}

    $$\begin{array}{|l|l|l|l|l|}

\multicolumn{5}{l}{\text{(a) Predicted formulas for the nontrivial
part of $\Nha(p)$}}\\
    \hline
        a \setminus h & \ANY & \PR & \RP & \RPPR  \\
        \hline
    \ANY & \approx
    \sum\frac{\abs{S_{m}}^{2}}{\abs{T_{m}}}  &
    \approx \scriptstyle \phi(p-1) &
    \approx \scriptstyle \phi(p-1) &
    \approx \frac{\phi(p-1)^{3}}{(p-1)^{2}} \\

        \hline
    \PR & \approx \scriptstyle \phi(p-1) &
    \approx \frac{\phi(p-1)^{2}}{(p-1)} &
    \approx \frac{\phi(p-1)^{2}}{(p-1)} &
    \approx \frac{\phi(p-1)^{3}}{(p-1)^{2}} \\

    \hline
    \RP &  \approx \scriptstyle \phi(p-1) &
    \approx \frac{\phi(p-1)^{2}}{(p-1)} &
    \approx \frac{\phi(p-1)^{2}}{(p-1)} &
    \approx \frac{\phi(p-1)^{3}}{(p-1)^{2}} \\

    \hline
    \RPPR & \approx \frac{\phi(p-1)^{3}}{(p-1)^{2}} &
    \approx \frac{\phi(p-1)^{3}}{(p-1)^{2}} &
    \approx \frac{\phi(p-1)^{3}}{(p-1)^{2}} &
    \approx \frac{\phi(p-1)^{3}}{(p-1)^{2}} \\

        \hline
\multicolumn{5}{l}{}\\
\multicolumn{5}{l}{\text{(b) Predicted values for the nontrivial
part of $\Nha(100057)$}}\\
        \hline
        a \setminus h & \ANY & \PR & \RP & \RPPR  \\
        \hline
\ANY & 190822.0 & 30240 & 30240 & 2762.225  \\
\hline
\PR & 30240 & 9139.458 & 9139.458 & 2762.225  \\
\hline
\RP & 30240 & 9139.458 & 9139.458 & 2762.225  \\
\hline
\RPPR & 2762.225 & 2762.225 & 2762.225 & 2762.225 \\
    \hline
\multicolumn{5}{l}{}\\
\multicolumn{5}{l}{\text{(c) Observed values for the nontrivial
part of $\Nha(100057)$}}\\
        \hline
        a \setminus h & \ANY & \PR & \RP & \RPPR  \\
        \hline
\ANY &  190526 & 30226 & 30291 & 2820 \\
\hline
\PR &   30226 & 9250 & 9231 & 2820 \\
\hline
\RP &   30291 & 9231 & 9086 & 2820\\
\hline
\RPPR & 2820 & 2820 & 2820 & 2820 \\
\hline
\end{array}$$
\end{table}

As observed in Section~\ref{equivsec}, to estimate the number of
solutions to~\eqref{tc} in the remaining cases we need to look
at~\eqref{hasmith}.  We start by estimating the number of nontrivial
solutions.  This requires a finer version of Heuristic~\ref{xtoxxheur}
which takes $d=\gcd(h,a,p-1)$ into account.

\begin{heuristic} \label{finextoxxheur}
    Fix $d$, $e$ such that $e$ divides $p-1$ and $d$ divides $e$.
    Then the map $x \mapsto x^{x/d} \mod{p}$  is a random map even when restricted
    to a specific order and greatest common divisor,
     in the sense that
    for all $p$, given $y \in \set{1, \ldots, p-1}$ such that $\ord_p y =
    f/\gcd(e,f)$,
    then
    \begin{multline*}
    \#\set{x \in \set{1, \ldots, p-1} \st x^{x/d} \equiv y \mod{p},
    \ \gcd(x,p-1)=e, \ \ord_p x = f } \\
    \approx \frac{ \#\set{z \in \set{1, \ldots, p-1} \st
    \gcd(z,p-1)=e, \ \ord_p z = f } }%
    {\#\set{w \in \set{1, \ldots, p-1} \st \ord_{p} w = \ord_{p} y}}.
    \end{multline*}
\end{heuristic}

Now we can approximate the number of nontrivial solutions of~\eqref{hasmith}
using a similar birthday paradox argument to that used
in~\cite{Holden02} for Conjecture~\ref{conj7}.

By Heuristic~\ref{finextoxxheur}, we see that the nontrivial part of
$\Nha'_{h \ANY, a \ANY}(p)$ is equal to:
\begin{multline*}
    \sum_{d \mid p-1} \#
    \set{(h,a) \st h \neq a, \ \text{\eqref{hasmith}~holds},
    \ \gcd(h,a,p-1)=d } \\
    \begin{aligned}
	& = \sum_{d \mid p-1} \sum_{\substack{e,f \mid p-1 \\ \gcd(e,f)=d}}
    \# \set{(h,a) \st h \neq a, \ \text{\eqref{hasmith}~holds},
    \ \gcd(h,p-1)=e, \ \gcd(a,p-1)=f } \\
    & \approx \sum_{d \mid p-1} \sum_{\substack{e,f \mid p-1 \\ \gcd(e,f)=d}}
    \sum_{m \mid p-1} \frac{\#S_{m, e} \cdot \#S_{m, f}}{\# T_{m} },
    \end{aligned}
\end{multline*}
 where
\begin{equation*}
    \begin{split}
    S_{m, r} &= \set{ x \st \ord_p (x^{x/d}) = m,\ \gcd(x, p-1)=r } \\
    & = \bigcup_{n \mid (p-1)/m} \set{x \st \ord_{p}(x) = nm,
    \ \gcd\left(\frac{x}{d},nm\right)=n, \ \gcd(x, p-1)=r} \\
    & = \bigcup_{\substack{n \mid (p-1)/m \\
    \gcd\left(\frac{r}{d},  nm\right)=n}}
    \set{x \st \ord_{p}(x) = nm,
    \ \gcd(x, p-1)=r}
    \end{split}
\end{equation*}
and
\[
T_{m} = \set{ x \st  \ord_{p} x = m }.
\]
Then, by Heuristic~\ref{indepheur}, we have:
\begin{equation*}
    \begin{split}
    \# S_{m, r}
    & \approx \sum_{\substack{n \mid (p-1)/m \\
    \gcd\left(\frac{r}{d},  nm\right)=n}}
    \frac{1}{p-1}
    \#\set{x \st \ord_{p}(x) = nm } \cdot \#\set{x \st \gcd(x, p-1)=r}
    \\
    & = \sum_{\substack{n \mid (p-1)/m \\
    \gcd\left(\frac{r}{d},  nm\right)=n}}
    \frac{1}{p-1} \phi(nm) \phi\left(\frac{p-1}{r} \right)
    \end{split}
\end{equation*}

Thus
\begin{multline*}
    \sum_{d \mid p-1} \#
    \set{(h,a) \st h \neq a, \ \text{\eqref{hasmith}~holds},
    \ \gcd(h,a,p-1)=d } \\
    \begin{aligned}
    & \approx \sum_{d \mid p-1} \sum_{\substack{e,f \mid p-1 \\ \gcd(e,f)=d}}
    \sum_{m \mid p-1} \frac{1}{\phi(m)}\left(
    \sum_{\substack{n \mid (p-1)/m \\
    \gcd\left(\frac{e}{d},  nm\right)=n}}
    \frac{1}{p-1} \phi(nm) \phi\left(\frac{p-1}{e} \right) \right)
    \\
    & \phantom{\approx \sum_{d \mid p-1} \sum_{\substack{e,f \mid p-1 \\ \gcd(e,f)=d}}
    \sum_{m \mid p-1} \frac{1}{\phi(m)}}
    \times \left( \sum_{\substack{t \mid (p-1)/m \\
    \gcd\left(\frac{f}{d},  tm\right)=t}}
    \frac{1}{p-1} \phi(tm) \phi\left(\frac{p-1}{f} \right) \right) \\
    & = \sum_{d \mid p-1} \sum_{\substack{e,f \mid p-1 \\ \gcd(e,f)=d}}
    \sum_{m \mid p-1}
    \sum_{\substack{n \mid (p-1)/m \\
    \gcd\left(\frac{e}{d},  nm\right)=n}}
    \sum_{\substack{t \mid (p-1)/m \\
    \gcd\left(\frac{f}{d},  tm\right)=t}}
    \frac{\phi(nm) \phi(tm)
    \phi\left(\frac{p-1}{e} \right)
    \phi\left(\frac{p-1}{f} \right)}{(p-1)^{2}\phi(m)}.
    \end{aligned}
\end{multline*}

\begin{proposition}
    For any $d$ dividing $q$,
    \begin{equation*}
    \sum_{\substack{e,f \mid q \\ \gcd(e,f)=d}}
    \sum_{m \mid q}
    \sum_{\substack{n \mid q/m \\
    \gcd\left(\frac{e}{d},  nm\right)=n}}
    \sum_{\substack{t \mid q/m \\
    \gcd\left(\frac{f}{d},  tm\right)=t}}
    \frac{\phi(nm) \phi(tm)
    \phi\left(\frac{q}{e} \right)
    \phi\left(\frac{q}{f} \right)}{\phi(m)}
    = q \ J_{2}\left( \frac{q}{d} \right),
    \end{equation*}
    where $J_{2}(r)$ is the Jordan function $J_{2}(r) =\sum_{s \mid r}
    s^{2}\mu\left(\frac{r}{s}\right)$.
\end{proposition}

\begin{proof}
    This can be verified directly when $q$ is a prime power; then use
    multiplicativity for the general case.
\end{proof}

It seems likely that a more combinatorial proof of this proposition
can be found.

Finally, we see that the nontrivial part of $\Nha'_{h \ANY,
a \ANY}(p)$ is approximately

\begin{multline*}
    \sum_{d \mid p-1} \sum_{\substack{e,f \mid p-1 \\ \gcd(e,f)=d}}
    \sum_{m \mid p-1}
    \sum_{\substack{n \mid (p-1)/m \\
    \gcd\left(\frac{e}{d},  nm\right)=n}}
    \sum_{\substack{t \mid (p-1)/m \\
    \gcd\left(\frac{f}{d},  tm\right)=t}}
    \frac{\phi(nm) \phi(tm)
    \phi\left(\frac{p-1}{e} \right)
    \phi\left(\frac{p-1}{f} \right)}{(p-1)^{2}\phi(m)} \\
    \begin{aligned}
    & = \sum_{d \mid p-1} \frac{1}{p-1} \ J_{2}\left( \frac{p-1}{d} \right) \\
    & = p-1.
    \end{aligned}
\end{multline*}

As we saw in Section~\ref{equivsec}, Heuristic~\ref{uvheur} implies 
that every nontrivial solution to~\eqref{hasmith} with $h \ANY$ 
and $a \ANY$ produces on
average one pair $(g,h)$ which is a nontrivial solution to~\eqref{tc}.
Thus the nontrivial part of $\Ntc_{g \ANY, h \ANY}(p)$ and
also the nontrivial part of $\Nha'_{h \ANY, a \ANY}(p)$ are
both approximately equal to $p-1$.

Similarly, we saw that Heuristic~\ref{uvheur} implies that every
solution to~\eqref{hasmith} with $h \ANY$ 
and $a \ANY$ produces on average
$\phi(p-1)/(p-1)$
pairs $(g,h)$ which
are solutions to~\eqref{tc} with $g \PR$.
Combining this with the previous argument, we see that the nontrivial
part of $\Ntc_{g \PR, h \ANY}(p)$ is approximately
\[
\sum_{d \mid p-1} \frac{J_{2}\left( \frac{p-1}{d} \right)}{p-1}
\frac{\phi\left(p-1\right)}{p-1} = \phi(p-1).
\]

These calculations justify the following conjectures:
\begin{conjecture} \label{conj8} \mbox{}
    \begin{enumerate}
    \item
$\Ntc_{g \PR, h \ANY}(p) \approx 2\phi(p-1)$.
\item
    $\Ntc_{g \ANY, h \ANY}(p) \approx
2(p-1)$.
\end{enumerate}
\end{conjecture}

These conjectures were made in~\cite{Holden02} on the basis of an
extension of the ``random map'' idea for $x \mapsto \log x$.  As was
explained there, however, it was not clear how to formulate the idea
as a heuristic that could be proved in a rigorous form.  The new
analysis explains the complications in the relationship between
$\Nha_{h \ANY, a \ANY}$ and $\Ntc_{g \ANY, h \ANY}$
encountered in~\cite{Holden02}.

Finally, Heuristic~\ref{uvheur} can be used to justify the
last set of conjectures from~\cite{Holden02a}:

\begin{conjecture} \label{conj9} \mbox{} 
\begin{enumerate}
    \item $\Ntc_{g \RP, h \bullet}(p) \approx
    \left[\phi(p-1)/(p-1)\right] \Ntc_{g \ANY, h
    \bullet}(p)$.  \item $\Ntc_{g \RPPR, h \bullet}(p)
    \approx \left[\phi(p-1)/(p-1)\right] \Ntc_{g \PR, h \bullet}(p)$.
\end{enumerate}
\end{conjecture}

The conjectures on~\eqref{tc} are summarized in
Table~\ref{tctalktable}, which appeared in~\cite{Holden02a}.  The
table also contains new data collected since~\cite{Holden02}.  The
data sets from Tables~\ref{fptalktable}, \ref{hatalktable},
and~\ref{tctalktable} were collected on a Beowulf cluster with 19
nodes, each consisting of 2 Pentium~III processors running at 1~Ghz.
The programming was done in C, using MPI, OpenMP, and OpenSSL
libraries.  The collection took 68 hours for all values of
$\Nfp(p)$, $\Ntc(p)$, and $\Nha(p)$, for five primes $p$ starting at
100000.

\begin{table}[!ht]
    \caption{Solutions to~\eqref{tc}}
    \label{tctalktable}

   $$\begin{array}{|l|l|l|l|l|}
\multicolumn{5}{l}{\text{(a) Predicted formulas for the nontrivial
part of $\Ntc(p)$}}\\
    \hline
        g \setminus h & \ANY & \PR & \RP & \RPPR  \\
        \hline
    \ANY & \approx \scriptstyle (p-1)    &
    \approx \frac{\phi(p-1)^{2}}{(p-1)} &
    \approx \scriptstyle \phi(p-1)  &
    \approx \frac{\phi(p-1)^{3}}{(p-1)^{2}} \\

        \hline
    \PR & \approx \scriptstyle \phi(p-1)   &
     \approx \frac{\phi(p-1)^{2}}{(p-1)} &
    \approx \frac{\phi(p-1)^{2}}{(p-1)} &
   \approx \frac{\phi(p-1)^{3}}{(p-1)^{2}} \\

    \hline
    \RP &  \approx \scriptstyle \phi(p-1) &
    \approx \frac{\phi(p-1)^{3}}{(p-1)^{2}} &
     \approx \frac{\phi(p-1)^{2}}{(p-1)} &
   \approx \frac{\phi(p-1)^{4}}{(p-1)^{3}} \\

    \hline
    \RPPR & \approx \frac{\phi(p-1)^{2}}{(p-1)}     &
    \approx \frac{\phi(p-1)^{3}}{(p-1)^{2}} &
    \approx \frac{\phi(p-1)^{3}}{(p-1)^{2}} &
    \approx \frac{\phi(p-1)^{4}}{(p-1)^{3}} \\

        \hline
\multicolumn{5}{l}{}\\
\multicolumn{5}{l}{\text{(b) Predicted values for the nontrivial
part of $\Ntc(100057)$}}\\
        \hline
        g \setminus h & \ANY & \PR & \RP & \RPPR  \\
        \hline
\ANY & 100056 & 9139.5 & 30240 & 2762.2  \\
\hline
\PR &30240 & 9139.5 & 9139.5 & 2762.2 \\
\hline
\RP & 30240 & 2762.2 & 9139.5 & 834.8   \\
\hline
\RPPR & 9139.5 & 2762.2 & 2762.2 & 834.8  \\
    \hline
\multicolumn{5}{l}{}\\
\multicolumn{5}{l}{\text{(c) Observed values for the nontrivial
part of $\Ntc(100057)$}}\\
        \hline
        g \setminus h & \ANY & \PR & \RP & \RPPR  \\
        \hline
\ANY & 100860 & 9231 & 30291 & 2820\\
\hline
\PR & 30850 & 9231 & 9231 & 2820\\
\hline
\RP &    30368 & 2882 & 9240 & 916\\
\hline
\RPPR & 9376 & 2882 & 2882 & 916\\

    \hline
\end{array}$$
\end{table}

 \section{Averages of the main terms}

 \label{aveconjsec}
 
 Thus far we have considered variants of Brizolis conjecture for a
 fixed finite field with $p$ elements. In the next two sections we
 consider average versions of these results and conjectures.
 The conjectures predict a main term; the results give a main
 term and an error term. The following sequence of lemmas gives
 the behavior of the main terms, on average. The only
 result from analytic number theory we need in order
 to prove these lemmas is the so-called Siegel-Walfisz theorem.
 As usual $\pi(x;d,a)$ denotes the number of primes $p\le x$ such
 that $p\equiv a$ modulo $d$, and ${\Li}(x)=\int_2^x{dt/\ln t}$
 denotes the logarithmic integral. 
 \begin{lemma}[{\cite[Satz 4.8.3]{Prachar}}]
 \label{walvis}
  Let $C>0$ be arbitrary. Then
  \[
 \pi(x;d,a)=\frac{{\Li}(x)}{\phi(d)}+O(xe^{-c_1\sqrt{\ln x}}),
 \]
 uniformly for $1\le d\le \ln^Cx$, $(a,d)=1,$ where the constants depend
 at most on $C.$
 \end{lemma}
 The following result for $k=1$ is well-known, see e.g.
 \cite{MoreeJNT, Stephens}. For arbitrary $k$ it
 was claimed by
 Esseen \cite{Esseen} (but only proved for $k=3$). We
 present a proof based on an idea of Carl Pomerance \cite{PomerancePC}.
 An analogue of this result for natural numbers was proved by
 Issai Schur in his Winter Semester lectures of 1923-24. He proved,
 for any complex number $s$, that
 \[
 \lim_{m\rightarrow \infty}\frac{1}{m}
 \sum_{n=1}^m\left(\frac{\phi(n)}{n}\right)^s
 =\prod_p\left(1+\frac{(1-1/p)^s-1}{p}\right).
 \]
 For an instructive discussion of this result see
 \cite[Chapter 4.2]{Kac}.
 \begin{lemma}
 \label{kinderlijk1}
 Let $k$ and $C$ be arbitrary real numbers with $C>0$. Then
 \[
 \sum_{p\le x}\left(\frac{\phi(p-1)}{p-1}\right)^k=A_k~{\Li}(x)
 +O_{C,k}\left(\frac{x}{\ln^C x}\right),
 \]
 where
 \[
 A_k=\prod_p\left(1+\frac{(1-1/p)^k-1}{ p-1}\right).
 \]
 \end{lemma}
 
\begin{proof}
    (The implicit constants in this proof depend at most on $C$ and $k$.)
 Let $g_k$ be the Dirichlet convolution of the M\"obius function
 and $(\phi(n)/n)^k$. Notice that $g_k$ is a multiplicative
 function and that $(\phi(n)/n)^k=\sum_{d|n}g_k(d)$. Using
 the latter identity we infer that
 \[
 \sum_{p\le x}\left(\frac{\phi(p-1)}{p-1} \right)^k=
 \sum_{p\le x}\sum_{d|p-1}{g_k(d)}=
 \sum_{d\le x}g_k(d)\pi(x;d,1).
 \]
 If $p$ is a prime, then
 clearly $g_k(p)=(1-1/p)^k-1$ and $g_k(p^r)=0$ for $r\ge 2$.
 For every $k$ there exist a constant $c_k$ such that $|g_k(p)|\le c_k/p$ for
 every prime $p$. Note that
 \begin{equation}
 \label{eerste} 
 |g_k(n)|\le \frac{c_k^{\omega(n)}|\mu(n)|}{n}\ll n^{-1+\epsilon},
 \end{equation}
 where $\omega(n)$ denotes the number of distinct prime divisors of $n$.
 Now write
 \begin{eqnarray*}
 \sum_{d\le x}{g_k(d)}\pi(x;d,1)&=&\sum_{d\le \ln^Bx}{g_k(d)}
 \pi(x;d,1)+\sum_{\ln^Bx<d\le x}{g_k(d)}\pi(x;d,1)\\
 &=&S_1+S_2,
 \end{eqnarray*}
 say, where $B>0$ is arbitrary for the moment.  In order to estimate
 $S_1,$ we invoke Lemma~\ref{walvis}.  This gives 
 \[
 S_1={\Li}(x)\sum_{d\le \ln^Bx}\frac{g_k(d)}{ \phi(d)}+
 O_C\left(\frac{x}{\ln^C x}\right).
 \]
 Now
 \[
 \sum_{d\le \ln^Bx}\frac{g_k(d)}{ \phi(d)}=\sum_{d=1}^{\infty}
 \frac{g_k(d)}{ \phi(d)}+O\left(\sum_{d>\ln^Bx} \frac{|g_k(d)|}{\phi(d)}\right).
 \]
 We have $d/\phi(d)=\prod_{p|d}(1-p^{-1})^{-1}\le \prod_{p\le
 d}(1-p^{-1})^{-1}\ll \ln d,$ using Mertens' formula. This together
 with the estimate~\eqref{eerste} shows that
 the sum $\sum_{d=1}^{\infty}g_k(d)/\phi(d)$ is
 absolutely convergent.
 Since, moreover, $g_k(d)/\phi(d)$ is multiplicative, we find using
 the Euler product identity that
 $\sum_{d=1}^{\infty}{g_k(d)/\phi(d)}=
 A_k$.
 Using~\eqref{eerste} we infer that
  \[
 \sum_{d>\ln^Bx}
 \frac{|g_k(d)|}{\phi(d)}\ll \sum_{d>\ln^Bx}\frac{\ln d}{d^{2-\epsilon}}
 \ll \frac{B\ln \ln x}{\ln^{B(1-\epsilon)} x}.
 \]
 Invoking the estimates
 $\pi(x;d,1)<x/d$, and~\eqref{eerste} leads to $S_2=O(x\ln^{-B(1-\epsilon)}x).$
 On putting everything together and taking $B$ sufficiently large, the result
 follows.
 \end{proof}

\begin{remark} \label{constantsrem}
    Using, e.g., Maple it turns out that in the range
 $0\le k\le 27$ the constant $A_k$ is quite well approximated
 by $e^{-1.011k+0.0278k^2}$. The constant $A_1$ equals the Artin constant. 
 Let
 \[
 A_{k,n}=\prod_{p>n}\left(1+\frac{(1-1/p)^k-1}{p-1}\right)\text{~and~}
 \zeta_n(k)=\zeta(k)\prod_{p\le n}(1-p^{-k}).
 \]
 If $k$ is a natural number and $n$ is sufficiently large, then
 $A_{k,n}=\prod_{k\ge 2}\zeta_n(r)^{e_{k,r}}$, where the exponents
 $e_{k,r}$ are integers that can be explicitly computed
 \cite{Moreeconstant}.  In this way $A_k$ and indeed any other Euler
 product appearing in this paper can be evaluated with arbitrary
 precision, cf.\  Theorem 2 of \cite{Moreeconstant}.  In
 Table~\ref{Aktable} we present a few examples.
 \end{remark}

 \begin{table}[!ht]
    \caption{The constants $A_k$}
    \label{Aktable}
\begin{center}
\begin{tabular}{|l|l|}\hline
$k$&$A_k$\\
\hline\hline
1&$0.37395~58136~19202~28805\cdots$\\ \hline
2&$0.14734~94000~02001~45807\cdots$\\ \hline
3&$0.06082~16551~20305~08600\cdots$\\ \hline
4&$0.02610~74463~14917~70808\cdots$\\ \hline
5&$0.01156~58420~47143~35542\cdots$\\ \hline
6&$0.00525~17580~26977~39754\cdots$\\ \hline
7&$0.00243~02267~63032~72703\cdots$\\ \hline
\end{tabular}
\end{center}
 \end{table}

If $a$ and $b$ are natural numbers, then by $(a,b)$ we denote the
greatest common divisor of $a$ and $b$ and by $[a,b]$ the lowest
common multiple.
\begin{lemma}
\label{viezig}
Let $a$ and $b$ be natural numbers and $C>0$. We have
\begin{equation}
\label{anderhalf}
 \sum_{\substack{p\le x \\ p\equiv 1 \mod{a} \\  p\equiv 1 \mod{b}}}
 \frac{\phi(\frac{p-1}{a})\phi(\frac{p-1}{b})}{(p-1)^2}
 =r(a,b)A_2{\Li}(x)
 +O_{a,b,C}\left(\frac{x}{\ln^C x}\right),
 \end{equation}
 where 
 \begin{equation}
 \label{derde}
r(a,b)=\frac{\phi(\frac{[a,b]}{(a,b)})\phi([a,b])}{[a,b]^4}
\prod_{p|ab}\frac{p(p^2+p-1)}{p^3-p^2-2p+1}
\prod_{p|\frac{ab}{(a,b)^2}}
\frac{p(p^2-1)}{p^3-2p+1}.
\end{equation}
We have
\begin{equation}
\label{vierde}
\frac{\phi(\frac{[a,b]}{(a,b)})}{\phi([a,b])[a,b]^2}
\le r(a,b)\le 4.13 \frac{\phi(\frac{[a,b]}{(a,b)})}{\phi([a,b])[a,b]^2}.
\end{equation}
\end{lemma}

\begin{proof}
    The proof can be carried out similarly to that of Lemma~\ref{kinderlijk1}.
We introduce an arithmetic function $h_{a,b}$ that satisfies
\begin{equation}
\label{tweede}
\frac{\phi(m \frac{a}{(a,b)})\phi(m\frac{b}{(a,b)})}%
{m\phi(\frac{a}{(a,b)})m\phi(\frac{b}{(a,b)})}
=\sum_{d|m}h_{a,b}(d).
\end{equation}
On noting that the left hand side of~\eqref{tweede} is a
multiplicative function of $m$, it follows that $h_{a,b}$ is
multiplicative.  Then $h_{a,b}$ is easily evaluated.  Taking
$m=(p-1)/[a,b]$ we find that~\eqref{anderhalf} holds with constant
\[
\frac{\phi(\frac{[a,b]}{(a,b)})}{\phi([a,b])[a,b]^2}
\sum_{d=1}^{\infty}\frac{h_{a,b}(d)\phi([a,b])}{\phi(d[a,b])}.
\]
After some manipulations the latter expression, in which the sum has
as argument a multiplicative function, is seen to equal
\begin{equation}
\label{Derde}
A_2 \frac{\phi(\frac{[a,b]}{(a,b)})}{\phi([a,b])[a,b]^2}
\prod_{p|ab} \frac{(p-1)(p^3-2p+1)}{p(p^3-p^2-2p+1)}
\prod_{p|\frac{ab}{(a,b)^2}} \frac{p(p^2-1)}{p^3-2p+1}.
\end{equation}
On further simplification this is seen to equal
$r(a,b)A_2$. It can be shown that
\[
\prod_p \frac{(p-1)(p^2-1)}{p^3-p^2-2p+1}\le 4.13.
\]
This inequality and the fact that the local factors in the two
products appearing in~\eqref{Derde} are all $>1$, then establishes the
truth of~\eqref{vierde}. 
\end{proof}

\begin{lemma}
\label{joshuaaverage}
Let $C>0$ be arbitrary. We have
\[
\sum_{p\le x} \frac{1}{p-1}
\sum_{e|p-1} \frac{1}{e}\phi\left(\frac{p-1}{e}\right)
=S~{\Li}(x)+O_C\left(\frac{x}{\ln^C x}\right),
\]
where 
\[
S=\prod_p\left(1-\frac{p}{p^3-1}\right)\approx 0.57595~99688~92945~43964\cdots
\]
is the Stephens constant (see~\cite{Stephens2}).
\end{lemma}
\begin{proof}
Using the fact that $\phi(n)/n=\sum_{d|n}\mu(d)/d$ with $n=(p-1)/e$,
we find on making the substitution $de=v$ and swapping the order of
summation that 
\[
\sum_{p\le x} \frac{1}{p-1}
\sum_{e|p-1} \frac{1}{e}\phi\left(\frac{p-1}{e}\right )
=\sum_{v\le x-1} \frac{\sum_{d|v}\mu(d)d}{v^2}
\sum_{\substack{p\le x \\ p\equiv 1 \mod{v}}}1.
\]
On splitting the summation range in the range $v\le \ln^B x$ and
$v>\ln^B x$ for an appropriate $B$, the result is then deduced as in
Lemma~\ref{kinderlijk1}.
\end{proof}

\begin{remark} Let $V=\{V_n\}_{n=0}^{\infty}$ be a sequence of
integers.  We say that $m$ divides the sequence $V$ if $m$ divides at
least one term of the sequence.  Denote by $\delta(V)$ the natural
density of primes $p$ dividing $V$, if it exists.  Stephens
\cite{Stephens2} proved, subject to the Generalized Riemann Hypothesis
(GRH), that $\delta(V)$ exists for a large class of second order
linear recurrences.  Moreover he showed, subject to GRH, that for
these sequences $\delta(V)$ equals a rational number times the
Stephens constant.  His work is extended and corrected in \cite{MS1,
MS2}.  For more details on the numerical approximation to $S$ given in
the lemma see \cite[p.~397]{Moreeconstant}.
\end{remark}

\begin{lemma}
\label{brizolisaverage}
Let $C>0$ be arbitrary. We have
\[
\sum_{p\le x}\frac{1}{(p-1)^2}
\sum_{e|p-1}\phi\left(\frac{p-1}{e}\right)^2
=A_1\frac{\zeta(3)}{\zeta(2)}{\Li}(x)
+O_C\left(\frac{x}{\ln^C x}\right),
\]
where 
\[
A_1 \frac{\zeta(3)}{\zeta(2)}
=\prod_p\left(1-\frac{2p}{p^3-1}\right)\approx 0.27327~30607~85299~15983\cdots.
\]
\end{lemma}
\begin{proof} Using the fact that $(\phi(n)/n)^2=\sum_{d|n}g_2(d)$
with $n=(p-1)/e$ (for the definition of $g_2(d)$ see the proof of
Lemma~\ref{kinderlijk1}), we find on making the substitution $de=v$
and swapping the order of summation that
\[
\sum_{p\le x} \frac{1}{p-1}
\sum_{e|p-1}\phi\left(\frac{p-1}{e}\right)^2
=\sum_{v\le x-1} \frac{\sum_{d|v} d^2 g_2(d)}{v^2}
\sum_{\substack{p\le x \\ p\equiv 1  \mod{v} }} 1.
\]
On splitting the summation range in the range $v\le \ln^B x$ and
$v>\ln^B x$ for an appropriate $B$, the result is then deduced as in
Lemma~\ref{kinderlijk1}.
\end{proof}

\begin{remark} Lemma \ref{viezig} suggests that the sum in the
previous lemma is asymptotically equal to
$A_2\sum_{e=1}^{\infty}r(e,e)$.  Some computation shows that, in
agreement with Lemma \ref{brizolisaverage}, we have
\[
A_2\sum_{e=1}^{\infty}r(e,e)
=A_2\sum_{e=1}^{\infty}\frac{1}{e^3}\prod_{p|e} \frac{p^3-2p+1}{p^3-p^2-2p+1}
=A_1 \frac{\zeta(3)}{\zeta(2)}.
\]
\end{remark}

\begin{lemma}
\label{U}
Let $C>0$ be arbitrary. We have
\begin{equation}
\label{upperbound}
\sum_{p\le x} \frac{1}{p-1}
\sum_{m|p-1}\phi(m)\left(\sum_{d|\frac{p-1}{m}}
\frac{\phi(d)}{d}\right)^2
=U~{\Li}(x)+O_C\left(\frac{x}{\ln^C x}\right),
\end{equation}
where 
\[
U=\prod_p\left(1+\frac{3p^2+2p+1}{p(p+1)(p^2-1)}\right)
\approx 3.4210\cdots.
\]
\end{lemma}

\begin{proof}
Let us define $h_1(n)=(\sum_{d|n}\phi(d)/d)^2$.  Note that $h_1$ is
multiplicative.  Let us denote the left hand side
of~\eqref{upperbound} by $I_1$.  We have
\begin{eqnarray*}
I_1&=&\sum_{p\le x}\sum_{m|p-1} \frac{\phi(\frac{p-1}{m})}{m\frac{p-1}{m}}
h_1(m) \\
&=& \sum_{p\le x}\sum_{m|p-1} \frac{h_1(m)}{m}
\sum_{\delta| \frac{p-1}{m}} \frac{\mu(\delta)}{\delta} \\
&=&\sum_{v\le x}\frac{\sum_{\delta|v}\mu(\delta)h_1(\frac{v}{\delta})}{v}
\sum_{\substack{p\le x \\ p\equiv 1 \mod{v}}} 1.
\end{eqnarray*}

Proceeding as in most of the earlier lemmas, we then deduce
that~\eqref{upperbound} holds true with constant
\[
\sum_{v=1}^{\infty}
\frac{\sum_{\delta|v}\mu(\delta)g(\frac{v}{\delta})}{v\phi(v)}
=\prod_p\left(1+\sum_{k=1}^{\infty}\frac{2+(2k-1)(1-1/p)}{p^{2k}}\right),
\]
which, after some tedious calculation, is seen to equal $U$. 
\end{proof}

\begin{lemma}
\label{L}
Let $C>0$ be arbitrary. We have
\begin{equation}
\label{lowerbound}
\sum_{p\le x} \frac{1}{p-1}
\sum_{m|p-1} \frac{\phi(m)^3}{m^2}\left(\sum_{d| \frac{p-1}{m}}
\frac{\phi(d)}{d}\right)^2
=L~{\Li}(x)+O_C\left(\frac{x}{\ln^C x}\right),
\end{equation}
where 
\[
L=\prod_p\left(1+\frac{p^5+2p^4-3p^3+p^2+1}{p^3(p+1)^3(p-1)}\right)
\approx 1.4446\cdots.
\]
\end{lemma}

\begin{proof} Let us denote the left hand side of~\eqref{lowerbound}
by $I_2$.  We have
\begin{eqnarray*}
I_2&=&
\sum_{p\le x}\sum_{m|p-1}
\frac{h_1(m)}{m}\left(\frac{\phi(\frac{p-1}{m})}{\frac{p-1}{m}}\right)^3 \\
&=& \sum_{p\le x}\sum_{m|p-1}\frac{h_1(m)}{m}
\sum_{\delta|\frac{p-1}{m}} g_3(\delta)\\
&=&\sum_{v\le x-1}
\frac{\sum_{\delta|v}\delta h_1(\delta)g_3(\frac{v}{\delta})}{v}
\sum_{\substack{p\le x \\ p\equiv 1 \mod{v}}} 1.
\end{eqnarray*}

Proceeding as in most of the earlier lemmas, we then deduce
that~\eqref{lowerbound} holds true with constant
\begin{multline*}
\sum_{v=1}^{\infty}
\frac{\sum_{\delta|v}\delta h_1(\delta)g_3(\frac{v}{\delta})}{v\phi(v)} \\
= \prod_p\left(1+\sum_{k=1}^{\infty} \frac{(1+k(1-\frac{1}{p}))^2
+((1-\frac{1}{p})^3-1)p(1+(k-1)(1-\frac{1}{p}))^2}{p^{2k}}\right),
\end{multline*}
which, after some tedious calculation, is seen to equal $L$. 
\end{proof}

The final lemma we will present is actually used in our error terms 
and not our main terms, but it is of the same character as the others 
in this section.

\begin{lemma}
\label{sigmalem}
Let $k$ and $C$ be arbitrary real numbers with $C>0$, $k > 0$. Then
\[ 
\sum_{p\le x} \frac{\sigma_{k}(p-1)}{(p-1)^{k}} = T_k {\Li}(x)
+O_{C,k}\left(\frac{x}{\ln^C x}\right),
\]
where 
\[
\sigma_{k}(n) = \sum_{d | n} d^{k}\]
and
\[
T_{k}=\prod_p\left(1+\frac{p}{(p-1)(p^{k+1}-1)}\right).
\]
\end{lemma}

\begin{proof}
    Using the fact that $\sigma_{k}(n)/n^{k}
    =\sum_{d|n} d^{k}/n^{k} = \sum_{d|n} 1/d^{k}$,
    we see that
    \[ 
    \sum_{p\le x} \frac{\sigma_{k}(p-1)}{(p-1)^{k}} =
    \sum_{p\le x}\sum_{d|p-1}{\frac{1}{d^{k}}}=
    \sum_{d\le x}\frac{1}{d^{k}}\pi(x;d,1).
    \]
    On splitting the summation range in the range $v\le \ln^B x$ and
    $v>\ln^B x$ for an appropriate $B$, the result is then deduced as
    in Lemma~\ref{kinderlijk1}.
\end{proof}

We have not yet computed the constants $T_{k}$ using the techniques 
described in Remark~\ref{constantsrem}, but a rough approximation 
using Maple gives the results shown in Table~\ref{Tktable}.

\begin{table}[!ht]
   \caption{The constants $T_k$}
   \label{Tktable}
\begin{center}
\begin{tabular}{|l|l|}\hline
$k$&$T_k$\\
\hline\hline
1&$2.20386\cdots$\\ \hline
2&$1.38098\cdots$\\ \hline
3&$1.15762\cdots$\\ \hline
4&$1.07163\cdots$\\ \hline
5&$1.03397\cdots$\\ \hline
6&$1.01646\cdots$\\ \hline
7&$1.00808\cdots$\\ \hline
\end{tabular}
\end{center}
\end{table}

\section{Averages of the conjectures and results}

\label{averessec}

Given the lemmas from the previous section it is trivial to establish
average versions of some of our results.  For example, we have:

\begin{theorem} Let $C>0$ be arbitrary.  We have
    \[
    \sum_{p\le x}{\frac{\Nfp_{g \PR,h \RPPR}(p)}{p-1}}
    =A_{2}{\Li}(x)+O_C\left(\frac{x}{\ln^C x}\right).
    \]
\end{theorem}

\begin{proof}
    Follows at once from Theorem~\ref{cz1thm}, Lemma~\ref{kinderlijk1},
    and the observation that, for every $\epsilon>0$, $\sum_{p\le
    x}\sigzero(p-1)\sqrt{p}(1 +\ln p)/(p-1)=O(x^{1/2+\epsilon})$. 
\end{proof}

Similarly, we have:

\begin{theorem} Let $C>0$ be arbitrary.  We have
\[
\sum_{p\le x}
\frac{\Mfp_{g{\PR},h{\ANY}}(p)}{p-1}
=A_1 \frac{\zeta(3)}{\zeta(2)}
{\Li}(x)+O_C\left(\frac{x}{\ln^C x}\right)
\]
and
\[
\sum_{p\le x}
\frac{\Mfp_{g{\ANY},h{\ANY}}(p)}{p-1}
=S~{\Li}(x)+O_C\left(\frac{x}{\ln^C x}\right).
\]
\end{theorem}

\begin{proof}
    Likewise follows from Theorems~\ref{thm5} and~\ref{thm6} and 
    Lemmas~\ref{joshuaaverage} and~\ref{brizolisaverage}. 
\end{proof}

Propositions~\ref{anyanyprop} and~\ref{pranyprop} are unfortunately
more problematic, due to the presence of the exceptionally large error
term.  As remarked there, the factor of $\sigone (p-1)-3(p-1)/2$ in the
error term can be averaged as
\begin{equation*}
\begin{split}
\sum_{p\le x} \frac{\sigone (p-1)-3(p-1)/2}{p-1} &=
(T_{2} - 3/2) \Li(x) +O_C\left( \frac{x}{\ln^C x}\right) \\
& \approx 0.70386
\Li(x) +O_C\left(\frac{x}{\ln^C x}\right).
\end{split}
\end{equation*}
(Apply Lemma~\ref{sigmalem}.)
The factor of $\sqrt{p}$, however, will still result in a error term 
with an order of magnitude larger than the main term.

On the other hand, almost all of the conjectures on~\eqref{fp},
\eqref{ha}, and~\eqref{tc} lend themselves easily to average versions
of the sort treated above.  For instance, we have:

\begin{conjecture} \mbox{}
    \begin{enumerate}
	\item     $\displaystyle \sum_{p\le x}{\frac{\Nfp_{g \ANY, h \ANY}(p)}{p-1}}
     \approx {\Li}(x).$
     
     	\item     $\displaystyle \sum_{p\le x}{\frac{\Nfp_{g \PR, h \ANY}(p)}{p-1}}
     \approx A_{1} {\Li}(x).$
     \end{enumerate}
\end{conjecture}

These conjectures and the average versions of our other conjectures
are summarized in Tables~\ref{fpavetable}, \ref{haavetable},
and~\ref{tcavetable}.  The data sets in these tables were collected on
the same Beowulf cluster with similar software.  The collection took
17 hours for all values of $\sum_{p\leq x}\frac{\Nfp(p)}{p-1}$,
$\sum_{p\leq x}\frac{\Ntc(p)}{p-1}$, and $\sum_{p\leq
x}\frac{\Nha(p)}{p-1}$, for $x=6143$.

The results of the preceding section unfortunately do not allow
us to evaluate the average value of the right hand side of 
Conjecture~\ref{conj7}(\ref{conj7a}).
Let us put
\[
w(p)=\sum_{m|p-1}\phi(m)\left(\sum_{d|m} \frac{\phi(dm)}{dm}\right)^2.
\]
Numerically it seems that
\[
\lim_{x\rightarrow \infty}\frac{1}{\pi(x)}
\sum_{p \leq x} \frac{w(p)}{p-1}
=1.644\cdots,
\]
with rather fast convergence. We are thus tempted to propose the
following conjecture.

\begin{conjecture}
 Let $C>0$ be arbitrary.  We have
  \[
 \sum_{p\leq x}
 \frac{\Nha_{a{\ANY},h{\ANY}}(p)}{p-1}
 =2.644\cdots{\Li}(x)+O_C\left(\frac{x}{\ln^C x}\right).
 \]
\end{conjecture}

Although we cannot prove (or even completely justify) this at present,
we can establish the following result.
\begin{lemma}
For every $x$ sufficiently large we have
\[
1.444\le \frac{1}{\pi(x)}
\sum_{p \leq x}\frac{w(p)}{p-1}\le 3.422
\]
\end{lemma}

\begin{proof}
Note that
\[
\sum_{m|p-1} \frac{\phi(m)^3}{m^2}\left(\sum_{d| \frac{p-1}{m}}
\frac{\phi(d)}{d}\right)^2\le
w(p)
\le \sum_{m|p-1}\phi(m)\left(\sum_{d|\frac{p-1}{m}}
\frac{\phi(d)}{d}\right)^2,
\]
where the first inequality, by the way, is exact if $p-1$ is squarefree.
The result now follows on invoking Lemma~\ref{L} and Lemma~\ref{U}.
\end{proof}

\begin{table}[!ht]
    \caption{Average Solutions to~\eqref{fp}}
    \label{fpavetable}
    $$\begin{array}{|l|l|l|l|l|}
\multicolumn{5}{l}{\text{(a) Predicted approximate values for $\frac{1}{\pi(x)}
\sum_{p \leq x} \Nfp(p)$}}\\
    \hline
        g \setminus h & \ANY & \PR & \RP & \RPPR  \\
        \hline
	\ANY & 1 & A_{2} & A_{1} & A_{2} \\

        \hline
    \PR & A_{1} & A_{2} & A_{2} & A_{2} \\

    \hline
    \RP & A_{1} & A_{3} & A_{2} & A_{3} \\

    \hline
    \RPPR & A_{2} & A_{3} & A_{3} & A_{3}\\
        \hline
\multicolumn{5}{l}{}\\
\multicolumn{5}{l}{\text{(b) Predicted approximate numeric values for $\frac{1}{\pi(x)}
\sum_{p \leq x} \Nfp(p)$}}\\
    \hline
        g \setminus h & \ANY & \PR & \RP & \RPPR  \\
        \hline
	\ANY & 1 & 0.1473494000 & 0.3739558136 & 0.1473494000 \\

        \hline
    \PR &0.3739558136 & 0.1473494000 &0.1473494000 & 0.1473494000 \\

    \hline
    \RP & 0.3739558136 & 0.0608216551 & 0.1473494000 & 0.0608216551 \\

    \hline
    \RPPR &0.1473494000 & 0.0608216551 & 0.0608216551 & 0.0608216551\\
        \hline
\multicolumn{5}{l}{}\\

    \multicolumn{5}{l}{\text{(c) Observed values for
    $x=6143$}}\\

        \hline
        g \setminus h & \ANY & \PR & \RP & \RPPR  \\
        \hline
   \ANY&  0.9904034375&0.14851987375&0.37592474125&0.14851987375\\

   \hline
\PR&   0.3749536975&0.14851987375&0.14851987375&0.14851987375\\
\hline
\RP&   0.3739629175&0.0612404775&0.15122619375&0.0612404775\\
\hline
\RPPR& 0.14792889125&0.0612404775&0.0612404775&0.0612404775\\

\hline
\end{array}$$
\end{table}

\begin{table}[!ht]
    \caption{Average Solutions to~\eqref{ha}}
    \label{haavetable}

    $$\begin{array}{|l|l|l|l|l|}

\multicolumn{5}{l}{\text{(a) Predicted approximate values for the nontrivial
part of}} \\
\multicolumn{5}{l}{\text{\mbox{}\qquad $\frac{1}{\pi(x)}
\sum_{p \leq x} \Nha(p)$}}\\  
    \hline
        a \setminus h & \ANY & \PR & \RP & \RPPR  \\
        \hline
    \ANY & 1.644\cdots & A_{1} & A_{1} &
    A_{3} \\

        \hline
    \PR &  A_{1} &  A_{2}& A_{2} &  A_{3} \\

    \hline
    \RP & A_{1} & A_{2} &  A_{2} &  A_{3} \\

    \hline
      \RPPR & A_{3} & A_{3} & A_{3}& A_{3}\\

    \hline
\multicolumn{5}{l}{}\\
\multicolumn{5}{l}{\text{(b) Predicted approximate numeric values for
the nontrivial part of}} \\
\multicolumn{5}{l}{\text{\mbox{}\qquad $\frac{1}{\pi(x)}
\sum_{p \leq x} \Nha(p)$}}\\  
    \hline
	a \setminus h & \ANY & \PR & \RP & \RPPR  \\
	\hline
    \ANY & 1.644\cdots &0.3739558136 & 0.3739558136 &
    0.0608216551 \\

	\hline
    \PR &  0.3739558136 & 0.1473494000& 0.1473494000 &  0.0608216551 \\

    \hline
    \RP & 0.3739558136 & 0.1473494000 &  0.1473494000 & 0.0608216551 \\

    \hline
      \RPPR & 0.0608216551 & 0.0608216551 & 0.0608216551& 0.0608216551\\

    \hline
\multicolumn{5}{l}{}\\
\multicolumn{5}{l}{\text{(c) Observed values for the nontrivial part 
for
    $x=6143$}}\\
	\hline
	a \setminus h & \ANY & \PR & \RP & \RPPR  \\
	\hline
\ANY & 1.6113896337 & 0.3655877485 & 0.3765792535 & 0.060552674 \\
\hline
\PR & 0.3655877485 & 0.14608992975 & 0.1478925015 & 0.060552674 \\
\hline
\RP & 0.3765792535 & 0.1478925015 & 0.146740421 & 0.060552674 \\
\hline
\RPPR &0.060552674 & 0.060552674 & 0.060552674 & 0.060552674\\
\hline
\end{array}$$
\end{table}

\begin{table}[!ht]
    \caption{Average Solutions to~\eqref{tc}}
    \label{tcavetable}

   $$\begin{array}{|l|l|l|l|l|}
\multicolumn{5}{l}{\text{(a) Predicted approximate values for the nontrivial
part of}} \\
\multicolumn{5}{l}{\text{\mbox{}\qquad $\frac{1}{\pi(x)}
\sum_{p \leq x} \Ntc(p)$}}\\  
    \hline
        g \setminus h & \ANY & \PR & \RP & \RPPR  \\
        \hline
    \ANY & 1 &  A_{2} & A_{1} & A_{3}\\

        \hline
    \PR & A_{1} & A_{2} & A_{2} & A_{3} \\

    \hline
    \RP & A_{1} & A_{3} & A_{2} & A_{4}\\

    \hline
    \RPPR & A_{2} & A_{3} & A_{3} & A_{4} \\

        \hline
\multicolumn{5}{l}{}\\

\multicolumn{5}{l}{\text{(b) Predicted approximate numeric values for the nontrivial
part of}} \\
\multicolumn{5}{l}{\text{\mbox{}\qquad $\frac{1}{\pi(x)}
\sum_{p \leq x} \Ntc(p)$}}\\  
    \hline
	g \setminus h & \ANY & \PR & \RP & \RPPR  \\
	\hline
    \ANY & 1 & 0.1473494000 & 0.3739558136 & 0.0608216551\\

	\hline
    \PR &0.3739558136 & 0.1473494000 &0.1473494000 & 0.0608216551 \\

    \hline
    \RP & 0.3739558136 & 0.0608216551 &0.1473494000 & 0.0261074463\\

    \hline
    \RPPR & 0.1473494000 & 0.0608216551 & 0.0608216551 & 0.0261074463 \\

	\hline
\multicolumn{5}{l}{}\\
\multicolumn{5}{l}{\text{(c) Observed values for the nontrivial part 
for
    $x=6143$}}\\
        \hline
        g \setminus h & \ANY & \PR & \RP & \RPPR  \\
        \hline
\ANY & 0.9933146575 & 0.14884923375 & 0.3772284725 & 0.06150940625\\
\hline
\PR & 0.37381320625 & 0.14884923375 & 0.14884923375 & 0.06150940625\\
\hline
\RP & 0.36701980375 & 0.06089004625 & 0.146029115 & 0.02640389625\\
\hline
\RPPR & 0.14697618875 & 0.06089004625 & 0.06089004625 & 0.02640389625\\

    \hline
\end{array}$$
\end{table}

\clearpage

\section{Conclusion and Future Work}

\label{conclsec}

Most of the theorems of Section~\ref{fpthmsec} suffer from an error
term which is larger than the main term.  This seems to be a direct
consequence of the use of Lemma~\ref{ordgcdlemma} and may be
unavoidable.  However, we have shown that we can put some limits on
how often the error actually approaches the worst case, and we have 
conjectured that even better limits exist.  The best next step  may be
further data collection in order to empirically count the number of primes 
with the potential for large errors.

We have begun to put our conjectures on a firm footing, deriving them
from as few heuristics as possible.  We hope to be able to prove these
heuristics in the future.  Then we should be able to convert the
conjectures into theorems by merely estimating the error term.

The project of extending our analysis to three-cycles and more
generally $k$-cycles for small values of $k$, mentioned
in~\cite{Holden02}, still remains to be done.  Along similar lines,
Igor Shparlinski has suggested attempting to analyze the average
length of a cycle, which could have many practical applications in the
analysis of cryptographically secure pseudorandom bit generators, as
mentioned in~\cite{Holden02}.

\section*{Acknowledgments}

Once again, the first author would like to thank the people mentioned
in~\cite{Holden02}: John Rickert, Igor Shparlinski, Mariana Campbell,
and Carl Pomerance.  He would also like to thank Victor Miller
for the suggestion to use the Smith Normal Form.  Both authors would 
like to thank the anonymous referees for many helpful comments.

\newcommand{\SortNoop}[1]{}

\end{document}